%% file: Vertical_Curves.tex
\documentclass[11pt,reqno]{amsart}
\usepackage[foot]{amsaddr}
\usepackage[utf8]{inputenc}
\usepackage[T1]{fontenc}
\usepackage{amsfonts}
\usepackage[all]{xy}
\usepackage[title,titletoc]{appendix}
\usepackage{graphicx}
\usepackage{amssymb}
\usepackage{amsthm}
\usepackage{mathtools}
\usepackage{comment}
\usepackage{lineno}
\usepackage{esint}
\usepackage{hyperref}
\usepackage{xcolor}
\hypersetup{
   colorlinks,
   linkcolor={red!80!black},
   citecolor={blue!70!black},
   urlcolor={blue!90!black}
}
\usepackage{microtype}
\usepackage{dutchcal}
\usepackage{upgreek}

\usepackage{tikz}
\usepackage{fourier}
\usepackage[mathscr]{euscript}
\usepackage[maxbibnames=99,backend=biber, sorting=nyt, doi=false, url=false]{biblatex}
\usepackage{pgfornament}

\newcommand{\HH}{\mathbb{H}}

\newcommand{\R}{\mathbb{R}}

\newcommand{\cH}{\mathcal{H}}
\newcommand{\ud}[0]{\,\mathrm{d}}

\makeatletter
\def\resetMathstrut@{%
  \setbox\z@\hbox{%
    \mathchardef\@tempa\mathcode`\(\relax
    \def\@tempb##1"##2##3{\the\textfont"##3\char"}%
    \expandafter\@tempb\meaning\@tempa \relax
  }%
  \ht\Mathstrutbox@1.2\ht\z@ \dp\Mathstrutbox@1.2\dp\z@
}
\makeatother
\newcommand{\0}{\mathbf{0}}

\newcommand{\from}{\colon}

\newtheorem{thm}{Theorem}[section]

\newtheorem{lemma}[thm]{Lemma}
\newtheorem{prop}[thm]{Proposition}
\newtheorem{cor}[thm]{Corollary}
\newtheorem{defn}[thm]{Definition}
\theoremstyle{remark}

\DeclareMathOperator{\VCone}{VCone}
\DeclareMathOperator{\Cone}{Cone}
\DeclareMathOperator{\Kor}{Kor}
\DeclareMathOperator{\diam}{diam}
\DeclareMathOperator{\vol}{vol}
\DeclareMathOperator{\id}{id}
\DeclareMathOperator{\supp}{supp}
\DeclareMathOperator{\inter}{inter}
\DeclareMathOperator{\closure}{clos}
\DeclareMathOperator{\Adj}{Adj}

\DeclareMathOperator{\diver}{div}

\title{Vertical curves and vertical fibers in the Heisenberg group}

\author{Gioacchino Antonelli}
\address[Gioacchino Antonelli]{Department of Mathematics, University of Notre Dame, Hurley Hall, 255 Hurley, Notre Dame, IN 46556, United States}
\email{gantonel@nd.edu}
\author{Robert Young}
\address[Robert Young]
{New York University, Courant Institute of Mathematical Sciences, 251 Mercer Street, 10012, New York, USA}
\email{ryoung@cims.nyu.edu}

\begin{document}
%\linenumbers
\input{intro.tex}

\input{verticalcurve.tex}

\input{contactdiffeo.tex}

\printbibliography[title={References}]

\end{document}

%% file: intro.tex
\begin{abstract}
    Let $\HH$ denote the three-dimensional Heisenberg group. In this paper, we study vertical curves in $\HH$ and fibers of maps $\HH \to \mathbb{R}^2$ from a metric perspective. We say that a set in $\HH$ is a vertical curve if it satisfies a cone condition with respect to a homogeneous cone with axis $\langle Z \rangle$, the center of $\HH$. This is analogous to the cone condition used to define intrinsic Lipschitz graphs.

    In the first part of the paper, we prove that connected vertical curves are locally biH\"older equivalent to intervals. We also show that the class of vertical curves coincides with the class of intersections of intrinsic Lipschitz graphs satisfying a transversality condition. Unlike intrinsic Lipschitz graphs, the Hausdorff dimension of a vertical curve can vary; we construct vertical curves with Hausdorff dimension either strictly larger or strictly smaller than 2. Consequently, there are intersections of intrinsic Lipschitz graphs with Hausdorff dimension either strictly larger or strictly smaller than 2.

    In the second part of the paper, we consider smooth functions $\beta$ from the unit ball $B$ in $\HH$ to $\mathbb{R}^2$. We show that, in contrast to the situation in Euclidean space, there are maps  such that $\beta$ is arbitrarily close to the projection $\pi$ from $\HH$ to the horizontal plane, but the average $\cH^2$ measure of a fiber of $\beta$ in $B$ is arbitrarily small.
\end{abstract}

\maketitle

\section{Introduction}

Over the past twenty years, a rich theory of rectifiability in the first Heisenberg group $\HH$ has been developed based on Lipschitz curves and intrinsic Lipschitz (iLip) graphs, see the seminal work \cite{FSSCILip, FSIlip}. By generalizations of Rademacher's theorem \cite{FSSCiLip2, VittoneiLip, AMiLip}, these submanifolds can be locally approximated by one-dimensional horizontal subgroups (Lipschitz curves) and two-dimensional vertical subgroups (iLip graphs). One gap in this theory, however, is the systematic study of non-horizontal curves; such curves occur, for instance, as fibers of maps from $\HH$ to $\R^2$ or as the intersections of iLip graphs, and they can often be locally approximated by the vertical subgroup $\langle Z\rangle$.
Some results on smooth vertical curves have been established in \cite{Kozhevnikov, KozhevnikovArticle, LeonardiMagnani, MagnaniStepanovTrevisan}, see Section \ref{sec:Previous} for a discussion. One of our main goals in this paper is to study new phenomena that arise in the non-smooth case.

In the first part of the paper, we study the class of $\lambda$--vertical curves.
A $\lambda$--vertical curve is a set that satisfies the following cone condition. We fix exponential coordinates $(x,y,z)$ on $\HH\equiv \mathbb R^3$. For $\lambda>0$, let
\begin{equation}\label{eqn:Vcone}
  \VCone_\lambda := \{(x,y,z)\in \HH : |z| \ge \lambda (x^2+y^2)\}.
\end{equation}
This is a homogenous cone centered on $\langle Z\rangle$; as $\lambda$ increases, the cone converges to $\langle Z\rangle$. A subset $E\subset \HH$ is a \emph{$\lambda$--vertical curve} if $E\subset p \VCone_\lambda$ for all $p\in E$. Notice that as $\lambda$ increases, $\VCone_\lambda$ becomes smaller and $\lambda$--vertical curves become closer to vertical lines. We should understand $\lambda$ as a measure of the \emph{verticality} of the curve.

It is not \emph{a priori} clear that a $\lambda$--vertical curve is a topological curve, but in Section~\ref{sec:Properties}, we show that any connected vertical curve is a topological curve, see Propositions~\ref{prop:intervals} and \ref{lem:biHolderonCompact}.
\begin{prop}\label{prop:Properties}
  If $E$ is a connected $\lambda$-vertical curve, then there are an interval $I\subset \R$ and a map $\gamma\from I\to \HH$ such that $\gamma$ sends $I$ homeomorphically to $E$ and 
  \begin{equation}\label{eq:monotone}
    \gamma(t)\in \gamma(s) \VCone^+_\lambda \qquad \text{for all $s,t\in I$ such that $s<t$,}
  \end{equation}
  where $\VCone^+_\lambda := \VCone_\lambda\cap \{z\geq 0\}$.

  If $E$ is compact, then we can take $\gamma$ to be a bi-H\"older map from $[0,1]$ to $E$. That is, there are exponents $0<\alpha<\beta<1$ depending on $\lambda$ and there is a $C>0$ depending on $E$ such that for all $s,t\in [0,1]$,
  $$C^{-1}|s-t|^\beta < d(\gamma(s),\gamma(t)) < C |s-t|^\alpha.$$
\end{prop}
If $\gamma\from I\to \HH$ is a map satisfying \eqref{eq:monotone}, we say that $\gamma$ is a \emph{$\lambda$--vertical map}.

The vertical curve condition is closely related to the cone condition used to define iLip graphs. We can define iLip graphs as follows. 
Let $d_{\Kor}$ be the Kor\'anyi distance on $\HH$, see \eqref{eq:koranyi}. Let $W\subset \HH$ be a two-dimensional vertical subgroup, i.e., $\langle Z\rangle \subset W$. For $\lambda \in (0,1)$, we let
\begin{equation}\label{eq:def-cone}
  \Cone_{W,\lambda} := \{p\in \HH : d_{\Kor}(p,W) \le \lambda d_{\Kor}(\0,p)\}.
\end{equation}

Then $\Cone_{W,\lambda}$ is a homogeneous cone centered on $W$, with angle depending on $\lambda$. Indeed, if $p=(x,y,0)$ is a nonzero horizontal vector, then $p\in \Cone_{W,\lambda}$ if and only if $|\sin \angle(p,W)| \le \lambda$.
%(see ???). 
For $\lambda \in (0,1)$ and $E\subset \HH$, we say that $E$ is a \emph{$\lambda$--iLip graph (over $W$)} if for all $p\in E$, $E\subset p \Cone_{W,\lambda}$. Notice that as $\lambda \to 1$, $\Cone_{W,\lambda}$ becomes larger. Here $\lambda$ plays the role of a Lipschitz constant: the bigger $\lambda$, the less restrictive the Lipschitz condition.

If $V=W^\perp$ is the horizontal orthogonal complement of $W$, then $V$ is a one-dimensional horizontal subgroup and $V\cap \Cone_{W,\lambda} = \{\0\}$, so each coset of $V$ intersects an iLip graph $E$ over $W$ at most once. That is, $E$ can be viewed as the graph of a function from a subset of $W$ to $V$. The group $\HH$ splits as $\HH=W\cdot V$; let $\Pi_V\from \HH \to V$ and $\Pi_W\from \HH\to W$ be the corresponding projections. 
If $\Pi_W(E)=W$, we say that $E$ is an \emph{entire} iLip graph. In this case, $\Pi_W$ sends $E$ homeomorphically to $W$.

We can characterize vertical curves in terms of iLip graphs in the following way. 
\begin{lemma}\label{lem:containment}
  The following two properties hold:
  \begin{enumerate}
  \item For every $\lambda>0$, there is an $L\in (0,1)$, and for every $L\in (0,1)$, there is a $\lambda>0$ such that if $E$ is a $\lambda$--vertical curve and $W$ is a vertical plane, then there is an entire $L$--iLip graph $\Gamma$ over $W$ such that $E\subset \Gamma$.
  \item For every $\lambda>0$, there is an $L\in (0,1)$, and for every $L\in (0,1)$, there is a $\lambda>0$ such that if $E\subset \HH$ and $E$ is an $L$--iLip graph over $W$ for every vertical plane $W$, then $E$ is a $\lambda$--vertical curve.
  \end{enumerate}
\end{lemma}

Furthermore, the theory of vertical curves lets us characterize certain intersections of intrinsic Lipschitz graphs.

\begin{thm}\label{thm:intersections}
  Let $W_1,W_2\subset \HH$ be two-dimensional vertical subgroups of $\HH$ such that $W_1\ne W_2$. If $\lambda>0$ is sufficiently large, then there are $L,L'\in (0,1)$ such that:
  \begin{enumerate}
    \item If $\Gamma_1$ and $\Gamma_2$ are entire $L$--iLip graphs over $W_1$ and $W_2$ respectively, then $\Gamma_{1}\cap \Gamma_{2}$ is a $\lambda$--vertical curve which is nonempty, closed, connected, and has no endpoints.
    \item If $E$ is a $\lambda$--vertical curve, then there are $\Gamma_1$ and $\Gamma_2$ that are entire $L'$--iLip graphs over $W_1$ and $W_2$ and $\Gamma_1\cap \Gamma_2$ is a vertical curve containing $E$. If $E$ is nonempty, closed, connected, and has no endpoints, then $\Gamma_1\cap \Gamma_2=E$.
  \end{enumerate}
\end{thm}

\begin{figure}[ht]
    \centering
\includegraphics[width=0.4\textwidth]{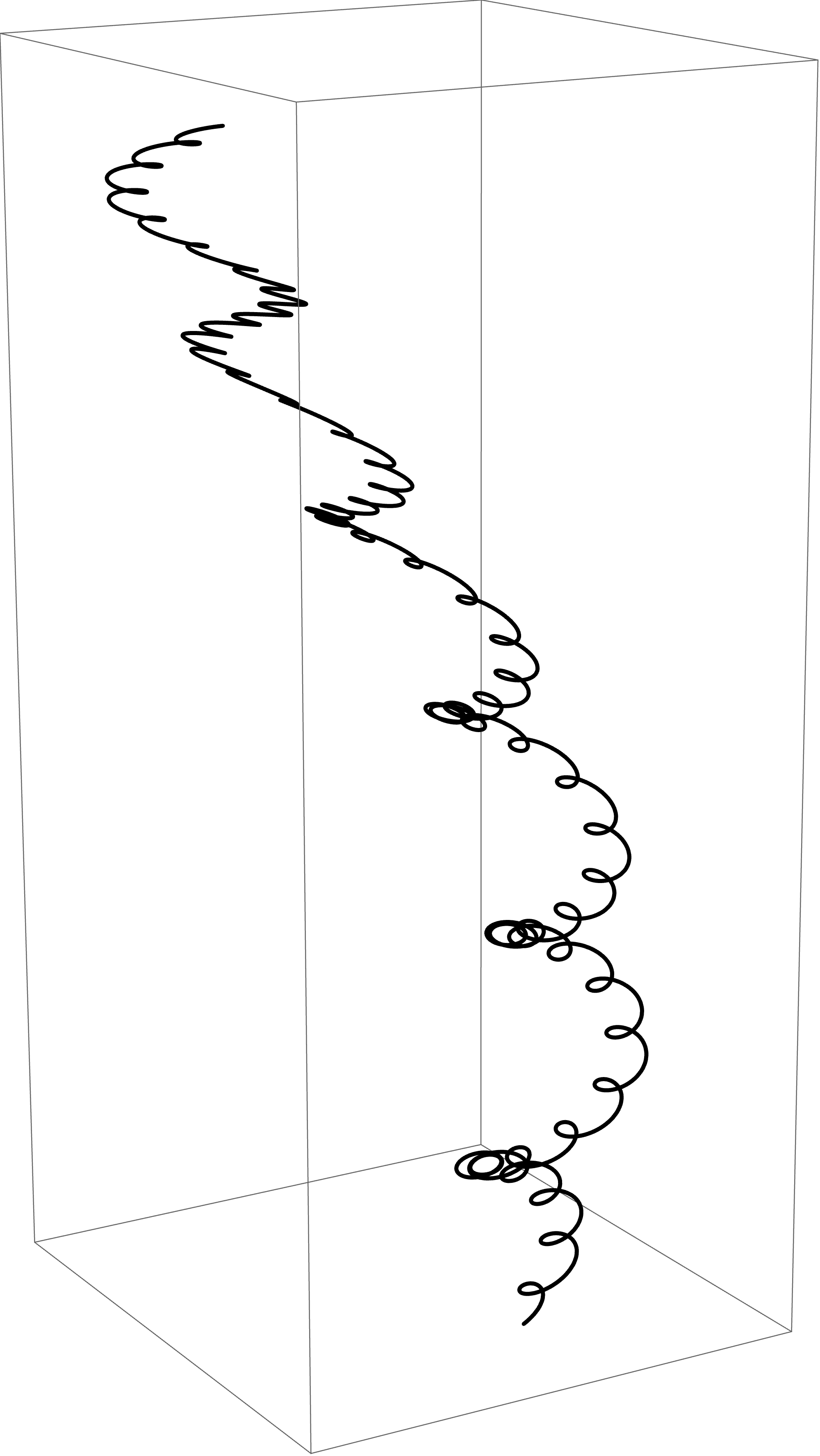}  
    \caption{An example of the curves constructed in Theorem \ref{thm:HausDim2}}
    \label{fig:Hausdim2}
\end{figure}

Finally, we construct $\lambda$--vertical curves with non-integer Hausdorff dimension. We construct these curves by starting with a smooth curve and applying helical perturbations at many different scales; an example can be seen in Figure~\ref{fig:Hausdim2}.
\begin{thm}\label{thm:HausDim2}
  For any $\lambda>0$, there are connected $\lambda$--vertical curves $E,E'\subset\HH$ with $\dim_H(E)<2<\dim_H(E')$.
\end{thm}
In fact, we show that if $\beta\from [0,1]\to \HH$ is a smooth curve such that $\beta'(t)$ is never a horizontal vector (and thus $\dim_H(\beta)=2$ as a direct consequence of Lemma~\ref{lem:measure-calculation}), then $\beta$ can be approximated arbitrarily closely by a vertical curve $\gamma$ with $\dim_H(\gamma)<2$ or $\dim_H(\gamma)>2$ (see Proposition~\ref{prop:curves}). 

The proof of Theorem~\ref{thm:intersections} 
 and Theorem~\ref{thm:HausDim2} immediately imply the following.
\begin{cor}
  For any two distinct vertical planes $W_1$ and $W_2$ and any $\lambda\in (0,1)$, there are entire $\lambda$--iLip graphs $\Gamma_1$ and $\Gamma_2$ over $W_1$ and $W_2$ such that $\dim_H(\Gamma_1\cap \Gamma_2)<2$ or $\dim_H(\Gamma_1\cap \Gamma_2)> 2$.
\end{cor}

These approximations and the curves constructed in Theorem~\ref{thm:HausDim2} are similar to Kozhevnikov's construction of vertical fibers connecting two points with arbitrarily large or arbitrarily small Hausdorff $2$--measure \cite{Kozhevnikov}.
A \emph{vertical fiber} is a curve that can be written in the form $f^{-1}(p)$ for some $C^1_H$ function $f\from \HH \to \R^2$ and some point $p\in \R^2$ which is a regular point of $f$. (Kozhevnikov called these ``vertical curves''; we have changed the name to avoid confusion.) Kozhevnikov proved that, in contrast to vertical curves, any vertical fiber is Reifenberg vanishing flat with respect to cosets of the center $\langle Z\rangle$ and thus has Hausdorff dimension $2$ \cite[Corollary 5.4.16]{Kozhevnikov}. Nonetheless, if $v,w\in \langle Z\rangle$, then there are vertical fibers that connect $v$ and $w$ and have arbitrarily large or small (possibly $\infty$ or $0$) $2$--dimensional Hausdorff measure \cite[Examples 5.6.16 \& 5.6.16]{Kozhevnikov}; see also the discussion in Section \ref{sec:Previous}. Kozhevnikov constructs these fibers using lacunary Fourier series. We will give a different construction for Theorem~\ref{thm:HausDim2}, but the geometric idea behind both constructions is that we can manipulate the Hausdorff $2$--measure of a curve by replacing smooth segments of the curve by helixes at many different scales.

\smallskip

In the second part of the paper, we study maps from $\HH$ to $\R^2$ and vertical fibers. Let $\pi \from \HH \to \R^2$ be the projection $\pi(x,y,z)=(x,y)$. 
We prove the following theorem.
\begin{thm}\label{thm:multiscale-maps}
  Let $B:=\overline{B}_1(\0)\subset \HH$ be the closed unit ball. 
  For any $\varepsilon>0$, there is a contact diffeomorphism $\beta \from B \to B$ with the following properties:
  \begin{enumerate}
  \item \label{it:identity} $\beta$ is the identity on $\partial B$. In fact, $\closure(\{x\in B:\beta(x)\neq x\})\subset \inter(B)$.
  \item \label{it:displace}
    For all $p\in \HH$, $d(p,\beta(p))<\varepsilon$.
  \item \label{it:small-fibers}
    \begin{equation}\label{eq:small-fibers}
      \int_{\R^2} \cH^2((\pi\circ \beta)^{-1}(v)) \ud v < \varepsilon.
    \end{equation}
  \end{enumerate}
\end{thm}
A contact diffeomorphism is a diffeomorphism that sends horizontal vectors to horizontal vectors; it is Lipschitz with respect to the Kor\'anyi metric. In particular, by Pansu's theorem, if $\alpha\from \HH\to \HH$ is a contact diffeomorphism, then for every $x\in \HH$, the differential $D\alpha_x$ is a Lie algebra isomorphism and thus $D[\pi\circ \alpha]_x$ is surjective on horizontal subspaces.

Therefore, if $\beta$ is as in Theorem~\ref{thm:multiscale-maps}, this implies that every $v\in \inter(\pi(B))$ is a regular point of $\pi\circ \beta$ and every curve $(\pi\circ \beta)^{-1}(v)$ is a vertical fiber. By \eqref{eq:small-fibers}, for all $v\in \pi(B)$ except a set of small measure, we have $\cH^2((\pi\circ \beta)^{-1}(v))\lesssim \varepsilon$.

Theorem~\ref{thm:multiscale-maps} contrasts with the behavior of maps from $\R^3$ to $\R^3$. Let $D^3\subset \R^3$ be the closed unit ball, $S^2=\partial D^3$, and let $\pi_{\R^2}\from \R^3\to \R^2$ be the orthogonal projection. Suppose that $f\from D^3\to D^3$ is a smooth map such that $f(x)=x$ for all $x\in \partial D$. Then for almost every $p\in D^2$, the preimage
$$\gamma_p := (\pi_{\R^2}\circ f)^{-1}(p)$$
contains a smooth curve connecting the two points of $\pi_{\R^2}^{-1}(p)\cap S^2$. Since $\pi_{\R^2}^{-1}(p) \cap D^3$ is the line segment between these two points, 
$$\cH^1(\gamma_p) \ge \cH^1(\pi_{\R^2}^{-1}(p)).$$
Therefore,
$$\int_{\R^2} \cH^1(\gamma_p) \ud p \ge \int_{\R^2} \cH^1(\pi_{\R^2}^{-1}(p)) \ud p = \vol(D^3).$$

\subsection{Methods}
Let us discuss the proofs of Theorem \ref{thm:HausDim2} and Theorem \ref{thm:multiscale-maps}. The construction of the sets $E,E'$ in Theorem \ref{thm:HausDim2} is based on the observation that any smooth curve $F\subset \HH$ whose tangent vectors are never horizontal is Reifenberg flat in the sense that $F$ is approximately vertical at small scales. That is, if $s_r\from \HH\to \HH$ is the automorphism that scales $\HH$ by a factor of $r$, then for any $p\in F$, the rescalings
$s_r(p^{-1}F)$ converge to $\langle Z \rangle$ as $r\to \infty$.

A helix around $\langle Z\rangle$ with the appropriate parameters is smooth, close to $\langle Z\rangle$, and has nonhorizontal tangents. Furthermore, depending on the direction of the helix, a segment of the helix is either longer or shorter than the corresponding segment of $\langle Z\rangle$. Since $F$ is close to vertical at small scales, we can perturb $F$ by replacing segments of $F$ by helixes. The result is a smooth curve with nonhorizontal tangents, so it is still Reifenberg flat. By repeating the process at smaller and smaller scales, we can make the $\cH^2$ measure as small or large as we want. By taking limits, we obtain vertical curves with Hausdorff dimension greater or less than $2$. Along the way, we can ensure that any curve constructed in this way is also a $\lambda$--vertical curve. 
\smallskip

The main idea of the proof of Theorem~\ref{thm:multiscale-maps} is the following, see Lemma~\ref{lem:reduction}. Given a contact map $\beta \from \HH\to\HH$ and a point $p\in \HH$ where the horizontal Jacobian of $\beta$ does not vanish, we show that one can perturb $\beta$ in a small  ball $B_{r_0}(p)$ so that the integral of the horizontal Jacobian on $B_{r_0}(p)$ decreases by a fixed factor $\alpha<1$. The perturbation can be taken to be a composition with a properly chosen contact diffeomorphism. By the Vitali Covering Theorem, we can apply perturbations on a collection of disjoint balls to reduce the integral of the horizontal Jacobian on $B$ by a factor $\alpha'<1$. Since the result is smooth, we can repeat the process at smaller scales to make the integral arbitrarily small. Theorem~\ref{thm:multiscale-maps} then follows from the co-area formula for smooth maps (Theorem \ref{thm:Coarea}).

\subsection{Structure}

 In Section \ref{sec:Properties} we prove Proposition \ref{prop:Properties}.
 In particular, Lemma~\ref{lem:Properties} and Proposition \ref{prop:intervals} prove several basic properties of $\lambda$--vertical curves, including the first part of Proposition \ref{prop:Properties}. We prove that vertical curves can be parametrized by maps that are bi-H\"older on compact sets in Lemma~\ref{lem:biHolderonCompact}, settling the second part of Proposition~\ref{prop:Properties}. In Section~\ref{sec:containment} we prove Lemma~\ref{lem:containment} and Theorem~\ref{thm:intersections}. In Section \ref{Proof2.5} we prove some auxiliary lemmas before proving Proposition~\ref{prop:curves}, from which Theorem~\ref{thm:HausDim2} follows.  In Section \ref{sec:FiberContact} we prove Lemma~\ref{lem:reduction} and use it to prove Theorem \ref{thm:multiscale-maps}.

 \subsection*{Acknowledgments} We would like to thank Davide Vittone for his fruitful discussions with us at the start of this project. G.A. has been partially supported by the AMS-Simons Travel Grant, and the NSF DMS Grant No. 2505713. R.Y. was supported by the National Science Foundation under Grant No.\ 2005609.

%% file: verticalcurve.tex
\section{Preliminaries and notation}

We recall here some preliminaries on the first Heisenberg group $\HH$ which will be useful throughout the paper. 

\subsection{The Heisenberg group and general notation}\label{sec:H1}

The first Heisenberg group $\HH$ is the unique simply connected Lie group whose Lie algebra is $\mathfrak{h}:=\mathrm{Lie}(\HH)=\langle X,Y\rangle \oplus \langle Z\rangle$, with $[X,Z]=[Y,Z]=0$ and $[X,Y]=Z$. We will identify $\HH$ with $\mathfrak{h}=\mathbb R^3$ by means of exponential coordinates. Namely, 
\[
\mathbb R^3 \ni (x,y,z) \leftrightarrow \exp(xX+yY+zZ) \in \HH.
\]
Under the above identification, $X,Y$, and $Z$ correspond to the elements $X=(1,0,0), Y=(0,1,0), Z=(0,0,1)\in \HH$.

We denote the coordinate functions on $\HH$ by $x,y,z\from \HH\to \R$. Let $\pi\from \HH \to \R^2$, $\pi(a) = (x(a),y(a))$ and let $\omega((x_1,y_1),(x_2,y_2)) = x_1y_2-y_1x_2$ be the area form on $\R^2$. Then $[u,v] = \omega(\pi(u),\pi(v))Z$.
We can write multiplication using the Baker--Campbell--Hausdorff formula; for $a,b\in \mathfrak{h}$, we have
$$\exp(a)\exp(b) = \exp\left(a+b+\frac{1}{2}[a,b]\right).$$
Using the identification of $\HH$ with $\mathfrak{h}$, we write
\begin{equation}
  \label{eq:group-law}
  a\cdot b = ab=a+b+\frac{1}{2}[a,b] = a + b + \frac{1}{2}\omega(\pi(a),\pi(b)) Z,
\end{equation}
or, in coordinates,
\[
  a b = (x_a,y_a,z_a)(x_b,y_b,z_b) := (x_a+x_b,y_a+y_b,z(ab)),
\]
where 
\[z(a b):=z_a+z_b+\frac{1}{2}\left(x_a y_b-x_b y_a\right).\]
(We use $ab$ and $a\cdot b$ interchangeably.)

The left-invariant vector fields on $\HH$ are generated by
\begin{align}\label{eq:left-inv}
  X & = \partial_x + \frac{y}{2} \partial_z & Y & = \partial_y - \frac{x}{2} \partial_z & Z & = \partial_z.
\end{align}
  By a slight abuse of notation, for $a\in \HH$ and $t\in \mathbb R$, we denote the one-parameter subgroup associated to $a$ by $\langle a \rangle$ and elements of $\langle a \rangle$ by
\[
a^t = (x_a,y_a,z_a)^t := \exp(t \log a) = (tx_a,t y_a, tz_a).
\]
This agrees with the usual exponent notation; in particular,
\[
a^{-1}=(x_a,y_a,z_a)^{-1}=(-x_a,-y_a,-z_a).
\]

Let us recall that the \textit{Kor\'anyi norm} of $(x,y,z)\in\HH$ is:
\begin{equation}\label{eq:def-koranyi}
  |(x,y,z)|_{\Kor}:=\sqrt[4]{(x^2+y^2)^2+16 z^2},
\end{equation}
and that we denote by
\begin{equation}\label{eq:koranyi}
d_{\Kor}(a,b):=|b^{-1}\cdot a|_{\Kor},
\end{equation}
the \textit{Kor\'anyi distance}. Note that for all $a\in \HH$, we have
\begin{equation}\label{eq:kor-proj}
  |\pi(a)| = \sqrt{x(a)^2+y(a)^2} \le |a|_{\Kor},
\end{equation}
with equality only if $z(a)=0$. We call the set $z^{-1}(0)\subset \HH$ the \emph{horizontal plane}.

When there is no risk of confusion we will just write $d$ instead of $d_{\Kor}$. The distance $d_{\Kor}$ is bi-Lipschitz equivalent to any Carnot--Carathéodory distance on $\HH$.

On $\HH$ we introduce a family of anisotropic scalings $\{S_r\}_{r>0}$, which are also group automorphisms, as follows:
\[
S_r(x,y,z):=(rx,ry,r^2z), \qquad \forall r>0, \quad \forall (x,y,z)\in \HH.
\]
These satisfy $d(S_r(a),S_r(b)) = r d(a,b)$ for all $r>0$ and all $a,b\in \HH$.

Finally, in this paper we will frequently use the notation
\[
A \lesssim B,
\]
to mean that there is a universal constant $C>0$ such that $A\leq CB$. Similarly for $A\gtrsim B$. When we write $A\approx B$ we mean that $B\lesssim A \lesssim B$. Moreover, when we write $A\lesssim_{\lambda} B$, we mean that the constant $C$ might depend on $\lambda$. Similarly for $\gtrsim_{\lambda}$, and $\approx_\lambda$.

\subsection{Intrinsic Lipschitz graphs in the Heisenberg group}

Let us recall some notation and definitions about intrinsic Lipschitz (iLip) graphs. As in \eqref{eq:def-cone}, if $W\subset \HH$ is a two-dimensional vertical subgroup, and $L\in (0,1)$, then we define
\[
\Cone_{W,L}:=\{p\in\HH:d_{\Kor}(p,W)\leq L d_{\Kor}(\0,p)\}.
\]
Then $\Cone_{W,L}$ is a scale-invariant cone centered on $W$.

An important example is the case of the $xz$--plane $W_0 = \{y=0\}$. This case is essentially general; rotation around the $z$--axis is an isometry of $\HH$, and any vertical plane is a rotation of the $xz$--plane. For any $p \in \HH$, we have $d_{\Kor}(p, W_0) = |y(p)|$; the distance from $W_0$ to $p$ is at least $|y(p)|$ by \eqref{eq:kor-proj}, and this bound is realized by 
$$d_{\Kor}(p, pY^{-y(p)}) = |Y^{-y(p)}|_{\Kor} = |y(p)|.$$
Then 
\[
\Cone_{W_0,L}=\{(x,y,z) \in\HH : |y| \le L\sqrt[4]{(x^2+y^2)^2+16z^2}\},
\]
i.e., $\Cone_{W_0,L}$ consists of points where $y$ is small compared to $x$ and $\sqrt{z}$. As $L\to 0$, $\Cone_{W_0,L}$ converges to $W_0$; as $L\to 1$, $\Cone_{W_0,L}$ gets bigger and bigger. Note, however, that if $t\ne 0$, then $Y^t\not \in \Cone_{W_0,L}$ for any $L\in (0,1)$.

The aperture of $\Cone_{W_0,L}$ depends on $L$; if $V(\theta) = (\cos \theta, \sin \theta, 0)$ is a horizontal unit vector, then $V(\theta) \in \Cone_{W_0,L}$ if and only if $|\sin \theta| \le L$, i.e., if the angle between $V(\theta)$ and $W_0$ is at most $\sin^{-1} L$.

The notion of iLip graph was introduced by Franchi--Serapioni--Serra Cassano in \cite{FSSCILip}. Here, following \cite{NY-VvsH} instead, we say that $E$ is an \emph{$L$--iLip graph (over $W$)} if for all $p\in E$, $E\subset p \Cone_{W,L}$. Moreover, we say that $E$ is an \emph{iLip graph (over $W$)} if there is $L\in (0,1)$ such that $E$ is an \emph{$L$--iLip graph (over $W$)}. There are many equivalent definitions of iLip graphs, but all of them give rise to the same class of objects. See, e.g., \cite[Section 3]{RigotiLip}. 

If $A\subset \HH$ is a one-parameter subgroup, we call $A$ \emph{horizontal} if $z(a)=0$ for all $a\in A$.
For every two-dimensional vertical subgroup $W\subset \HH$, we denote with $W^\perp\subset \HH$ the horizontal subgroup orthogonal to $W$. As we calculated for the $xz$--plane above, $W^\perp \cap \Cone_{W,L} = \{\0\}$ for any $L\in (0,1)$.

For each vertical subgroup $W\subset \HH$, $\HH$ splits as a product $\HH=W\cdot W^\perp$. Let $\Pi_W\from \HH\to W$ be the corresponding projection onto $W$. Then the fibers of $\Pi_W$ are the cosets of $W^\perp$. If $E$ is an iLip graph and $p\in E$, then
$$p\in E\cap p W^\perp \subset p (\Cone_{W,L}\cap W^\perp) =\{p\},$$
so each coset $p W^\perp$ intersects $p$ at most once. That is, $\Pi_W$ is injective on $E$.

We fix an orientation of $W^\perp$ by choosing a generator $\theta \in W^\perp$ of unit length. 
If $E$ is an iLip graph over $W$ and $U=\Pi_W(E)$, then there exists a function $\varphi\from U\to \R$ such that 
\[
E=\{w\cdot\theta^{\varphi(w)} : w\in U\}.
\]
We define
\[
E^+=\{w\cdot \theta^t:w\in U, t\geq \varphi(w)\}, \qquad E^-=\{w\cdot \theta^t:w\in U, t\leq \varphi(w)\}.
\]

We say that $E$ is an \emph{entire $L$-iLip graph (over $W$)} if it is an $L$-iLip graph (over $W$) and $\Pi_W(E)=W$. We will need the following extension theorem, which was proved in \cite{NY-VvsH}.
\begin{thm}[{\cite[Theorem 27]{NY-VvsH}}]\label{thm:extension}
  Fix $L\in (0,1)$ and a vertical plane $W\subset \HH$. Suppose $\Gamma$ is an intrinsic $L$--Lipschitz graph over $W$. Then there exists an entire intrinsic $L$--Lipschitz graph over $W$, denoted $\widetilde{\Gamma}$, such that $\Gamma\subset \widetilde{\Gamma}$.
\end{thm}
(In \cite{NY-VvsH}, this theorem is stated in terms of an extension of an intrinsic Lipschitz function $f\from U\to \R$, where $U\subset W$, but any iLip graph over $W$ can be written in this form by the discussion above.)

\subsection{Contact diffeomorphisms}\label{sec:Contact}
Let us recall some preliminaries on contact diffeomorphisms on $\HH$, which we will use in Section~\ref{sec:FiberContact}. Let $X,Y$, and $Z$ be the left-invariant vector fields on $\HH$ given by \eqref{eq:left-inv}. There is a contact structure on $\HH$ induced by the \emph{horizontal distribution}. This distribution associates to each point $v\in \HH$ the tangent plane $\langle X_v, Y_v\rangle \subset T_v\HH$.

For an open set $U\subset \HH$, we say that a smooth map $\beta\from U\to \HH$ is a \emph{contact map} if the differential of $\beta$ sends the horizontal distribution to the horizontal distribution, i.e., $D_v\beta(\langle X_v,Y_v\rangle) \subset \langle X_{\beta(v)},Y_{\beta(v)}\rangle$ for all $v\in U$.
A smooth vector field on $\HH$ is a \emph{contact vector field} if it generates a flow of \emph{contact maps}. By \cite{KoranyiReimann}, a vector field $V$ on $\HH$ is a contact vector field if and only if there exists $\psi\in C^\infty(\HH)$ such that
  $$
  V= V_\psi := Y[\psi] X - X[\psi] Y + \psi Z.
  $$
  (Here, $Y[\psi]=Y\psi$ is the directional derivative of $\psi$; we use square brackets to denote applying a differential operator when it makes an expression clearer.)
 For $\psi\in C^\infty(\HH)$, let us denote with $\Psi^t_\psi$ the flow of $V_\psi$ at time $t\in\mathbb R$. 

  Let us recall a useful computation. Let $\diver$ denote the Euclidean divergence. Then $X$, $Y$, and $Z$ are divergence-free, and by the product rule $\diver [f W] = f \diver W + W[f]$,
  $$
  \diver V_\psi = X[Y\psi] - Y[X\psi] + Z[\psi] = 2Z[\psi],
  $$
  where in the last equality we are using that $XY-YX=[X,Y]=Z$.
  Since the Euclidean volume form on $\R^3$ is equal to the Riemannian volume form on $\HH$, we can approximate the Jacobian determinant of $\Psi_\psi^t$ by using Liouville's formula:
  \begin{equation}\label{eqn:JacobianDet}
  J_{\Psi_\psi^t}(v) = 1+t \diver[V_\psi](v)+O(t^2)= 1 + 2Z[\psi](v) t + O(t^2),
  \end{equation}
  for all $v\in \HH$.

\section{Geometry and topology of vertical curves}\label{sec:Properties}
In this section, we prove some properties of vertical curves, including Proposition~\ref{prop:Properties}. Unlike intrinsic Lipschitz graphs, vertical curves cannot be naturally viewed as graphs over a subgroup of $\HH$. Indeed, if $\gamma\from I\to \HH$ is a vertical map, there need not be a $p$ such that $p\cdot\gamma(t)$ has nondecreasing $z$--coordinate.\footnote{Figure~\ref{fig:Hausdim2} shows one such curve. The key point is that if $\lambda(t)$ is a helix with large slope centered on the $z$--axis, then $\lambda(t)$ is a vertical curve, but $p\cdot \lambda(t)$ only has monotone $z$--coordinate when $p$ is near the $z$--axis. If $\gamma$ contains segments of two such helices, spaced far apart, then there is no $p$ such that $p\cdot \gamma(t)$ has monotone $z$--coordinate.}

Nonetheless, we will show that connected vertical curves are homeomorphic to intervals and that the homeomorphism can be taken to be bi-H\"older on compact sets.

We will need the following relation on $\HH$. For $p,q\in \HH$,
\begin{equation}\label{eq:def-order}
  p\prec q \qquad \text{if and only if } z(p^{-1}q)>0.
\end{equation}
Note that $z(p^{-1}q) = -z(q^{-1}p)$, so $p\prec q$ and $q \prec p$ are mutually exclusive. This is not a partial order on $\HH$ because it is nontransitive; one can show that $X \prec Y^{-1} \prec X^{-1} \prec Y \prec X$. Regardless, we will show that it induces a total order on connected vertical curves and use this to describe the geometry and topology of connected vertical curves.

The main results of this section are the following lemma and proposition.
\begin{lemma}\label{lem:Properties}
 Let $\lambda>0$, and let $E\subset \HH$ be a $\lambda$--vertical curve. Then:
  \begin{enumerate}
  \item For every $p\in\HH$, $p\cdot E$ is a $\lambda$--vertical curve.
  \item For every $\rho>0$, $s_\rho(E)$ is a $\lambda$--vertical curve.
  \item\label{Item3} If $\lambda > \frac{1}{4}$, or if $E$ is connected, then \eqref{eq:def-order} defines a total order on $E$. Moreover, if $E$ is the image of an injective map $\gamma\from I\to\HH$, then $\gamma$ is either monotone increasing in the sense that
  \begin{equation}\label{eqn:First}
  \gamma(s) \prec \gamma(t) \quad \forall s<t,
  \end{equation}
  or it is monotone decreasing in the sense that
  \begin{equation}\label{eqn:Second}
  \gamma(s) \succ \gamma(t) \quad \forall s<t.
  \end{equation}
  \item\label{Item4} 
    There is a $c_\lambda>0$ such that for all $a,b,c\in E$ such that $a\prec b$, $b\prec c$, and $a\prec c$,
    \begin{equation}\label{eq:coarse-monotone}
      z(a^{-1}\cdot c)\geq c_\lambda z(a^{-1}\cdot b) \text{ and } d(a,c)\geq c_\lambda d(a,b).
    \end{equation}
\end{enumerate}
\end{lemma}

\begin{prop}\label{prop:intervals}
  Let $\lambda>0$ and let $E$ be a connected $\lambda$--vertical curve. Then $E$ is homeomorphic to an interval and thus $E$ is parametrized by an increasing $\lambda$--vertical map $\gamma\from I \to\HH$.
\end{prop}
In fact, $E$ is locally bi-H\"older equivalent to an interval; see Lemma~\ref{lem:biHolderonCompact}.

We will first prove Lemma~\ref{lem:Properties}. Let
$$
\VCone_\lambda^+:=\VCone_\lambda\cap\{z\geq 0\}, \quad \VCone_\lambda^-:=\VCone_\lambda\cap\{z\leq 0\}.
$$ 

\begin{proof}[Proof of Lemma~\ref{lem:Properties}]
    Parts (1) and (2) are a direct consequence of the definitions, and the fact that $\VCone_\lambda$ is $s_\rho$-invariant.

Let us now prove part (3). Note that for every $p, q\in E$, we have \[q\in p\VCone_\lambda^+\cup p\VCone_\lambda^-,\] so $p\preceq q$ or $q\preceq p$. 
In order to prove part (3), we need to prove that the relation \eqref{eq:def-order} is transitive on $E$.

We first consider the case that $\lambda>\frac{1}{4}$. Let $a, b, c\in E$ and suppose $a\prec b$ and $b\prec c$. Then $g = a^{-1} b\in \VCone_\lambda^+$ and $h = b^{-1} c\in \VCone_\lambda^+$, i.e., $z(g) > \frac{1}{4}|\pi(g)|^2$ and $z(h) > \frac{1}{4}|\pi(h)|^2$. We claim that $a\prec c$, i.e., that $z(a^{-1}c) = z(gh) > 0$.

We have
\[
\begin{aligned}
z(gh) &= z(g) + z(h) + \frac{1}{2}\big(x(g) y(h) - x(h) y(g)\big) \\ &\ge z(g) + z(h) - \frac{1}{2}|\pi(g)||\pi(h)|,
\end{aligned}
\]
and
$$
\frac{1}{2}|\pi(g)||\pi(h)|<2\sqrt{z(g)z(h)}\le z(g)+z(h),
$$
by the arithmetic-geometric mean inequality. Therefore,
$$
z(gh) \ge z(g) + z(h) - \frac{1}{2}|\pi(g)||\pi(h)| > 0,
$$
which implies $a\prec c$.

Next, suppose that $E$ is connected. We again let $a, b, c\in E$ such that $a\prec b$ and $b\prec c$.
For $p\in E$, let 
\[
E_p^+:=\{x\in E : x\succ p\}, \quad E_p^-:=\{x\in E : x\prec p\}.
\]
Then $E\setminus\{p\}=E_p^+\sqcup E_p^-$, i.e., $E_p^+,E_p^-$ are disjoint, relatively open sets that cover $E\setminus\{p\}$. It follows that any connected component of $E\setminus\{p\}$ lies in one or the other. In particular, $a$ and $c$ are in different connected components of $E\setminus \{b\}$. 

Let $A$ be the connected component of $a$ in $E\setminus\{b\}$. Since $E$ is connected, the closure $\overline{A}$ of $A$ in $E$ contains $b$; since $a$ and $c$ are in different connected components of $E\setminus \{b\}$, $c\not \in \overline{A}$. Therefore, $\overline{A}$ is a connected subset of $E\setminus\{c\}$ that contains $a$ and $b$. Since $b\prec c$, this implies $a\prec c$, as desired.

Finally, let $\gamma\from I\to \HH$ be a $\lambda$--vertical map. Let $f\from I^2 \to \R$,  $f(s,t)=z(\gamma(s)^{-1}\gamma(t))$. Then $f(s,t)\ne 0$ for all $(s,t)$ such that $s\ne t$. In particular, either $f(s,t)>0$ for all $s<t$ or $f(s,t)<0$ for all $s<t$. In the first case, $\gamma(s)\prec \gamma(t)$ for all $s<t$; in the second case, $\gamma(s)\succ \gamma(t)$ for all $s<t$.

We now consider part (4). After a left-translation we may suppose that $a=\0$. Let $b=(x_b,y_b,z_b)$, and $c = (x_c,y_c,z_c)$. Then on one hand, since $c\succ b$, 
$$0\le z(b^{-1}c) = -z_b + z_c + \frac{1}{2}(-x_by_c+x_cy_b) \le -z_b + z_c + \frac{1}{2} |\pi(b)||\pi(c)|,$$
i.e.,  
\begin{equation}\label{eq-zc-lower-1}
  z_c \ge z_b - \frac{1}{2} |\pi(b)||\pi(c)|.
\end{equation}
Since $b\succ \0$, we have $|\pi(b)|\le \sqrt{z_b \lambda^{-1}}$. Thus, on one hand, if 
$|\pi(c)| \le \sqrt{z_b \lambda}$, then 
$$z_c \ge z_b - \frac{1}{2} \sqrt{z_b \lambda^{-1}}\sqrt{z_b \lambda} \ge \frac{z_b}{2}.$$
On the other hand, if $|\pi(c)| \ge \sqrt{z_b \lambda}$, then since $c\succ \0$, we have
$$z_c \ge \lambda(\sqrt{z_b \lambda})^2 = \lambda^2 z_b.$$
In either case, $z_c\ge \min\{\frac{1}{2},\lambda^2\} z_b$. 

Since $E$ is $\lambda$--vertical, if $p,q\in E$, then $d(p,q)\approx_\lambda |z(p^{-1}q)|$, i.e., $z_b\approx_\lambda d(a,b)$ and $z_c\approx_\lambda d(a,c)$. Therefore, there is a $c_\lambda>0$ such that $d(a,c)\ge c_\lambda d(a,b)$, as desired.
\end{proof}

The proof of  Proposition~\ref{prop:intervals} uses the following two lemmas.
\begin{lemma}\label{lem:compact-intersections}
  Let $\lambda>0$. There is a $C>0$ such that for all $p,q\in \HH$, if $q\in p\VCone_\lambda^+$, then
  $$\diam\left(p\VCone_\lambda^+\cap q \VCone_\lambda^-\right)\le C d(p,q).$$
  In particular, if $E$ is a $\lambda$--vertical curve and $p,v,q\in E$ are such that $p\prec v$, $v\prec q$, and $p\prec q$, then $d(p,v)\le C d(p,q)$.
\end{lemma}
\begin{proof}
  It suffices to prove the lemma when $p=\0$.
  Let $c_\lambda$ be as in Lemma~\ref{lem:Properties}.(\ref{Item4}), let $D=\VCone_\lambda^+\cap q \VCone_\lambda^-$, and let $v\in D$. Then $\0 \prec v$, $v\prec q$, and $\0 \prec q$, so by \eqref{eq:coarse-monotone}, $z(q) \ge c_\lambda z(v)$. By \eqref{eq:def-koranyi},
  $$
  z(v) \le c_\lambda^{-1} z(q) \le \frac{c_\lambda^{-1} d(\0,q)^2}{4}.
  $$
  Furthermore, since $v\in \VCone_\lambda^+$,
  $$
  |\pi(v)| \le \sqrt{z(v)\lambda^{-1}} \le d(\0,q) \sqrt{\frac{c_\lambda^{-1}\lambda^{-1}}{4}}.
  $$
  Therefore, if $C=\sqrt{c_\lambda^{-1}\lambda^{-1}/4} + \sqrt{c_\lambda^{-1}}$, then
  $$
  d(\0,v) = \sqrt[4]{|\pi(v)|^4+16 z(v)^2} \le C d(\0,q),
  $$
  as desired. 
\end{proof}

\begin{lemma}\label{lem:topologies}
  Let $\lambda>0$. Suppose that $E\subset \HH$ is a closed $\lambda$--vertical curve which is totally ordered, i.e., if $p,q,r\in E$, $p\prec q$, and $q\prec r$, then $p\prec r$. Then the
  order on $E$ is complete and the
  subset topology on $E$ coincides with the order topology, i.e., the topology generated by the rays
$I_{a,\infty}:=\{z\in E:a\prec z\}$ and $I_{-\infty,a}:=\{z\in E:z\prec a\}$, and the intervals $I_{a,b}:=\{z\in E:a\prec z\prec b\}$, with $a,b\in E$.
\end{lemma}
\begin{proof}
  (Unless we specify the order topology, topological notions in this proof will refer to the subset topology on $E$.)
  
  We first show that the order is complete, i.e., that every non-empty subset of $E$ with an upper bound has a least upper bound and every non-empty subset of $E$ with a lower bound has a greatest lower bound.

  Let $S$ be a subset of $E$ that is bounded above by some $M\in E$; the case that $S$ is bounded below is similar. Since $\HH$ is separable, $S$ contains a countable dense subset (in the subset topology), which we denote $\{a_1,\ldots,a_n,\ldots\}$.  For $i\ge 1$, let  $m_i=\max\{a_1,\ldots,a_i\}$. Then 
  $$
  m_i\in m_1\VCone_\lambda^+\cap M\VCone_\lambda^-, \quad \forall i\geq 0.
  $$ 
  By Lemma~\ref{lem:compact-intersections}, this intersection is compact, so there is a subsequence of $m_i$ that converges to some $v\in \HH$. Since $E$ is closed, $v\in E$, and one can check that $v$ is the least upper bound of $S$. Thus the order on $E$ is complete.

  Note that this argument also proves that if $S\subset E$ is bounded above, then $\sup S$ is a limit point of $S$ in the subset topology. 

  Now, we compare the order topology and the subset topology. For every $a,b\in E$ with $a\prec b$,
  $$I_{a,b} = E \cap \left(\HH \setminus (a \VCone_\lambda^- \cup b \VCone_\lambda^+)\right),$$
  so $I_{a,b}$ is open in the subspace topology (as are $I_{-\infty,a}$ and $I_{a,\infty}$). That is, the subspace topology is finer than the order topology.

  Conversely, let $p\in E$ and $\delta>0$. 
  We will show that there is an open set $U\subset E$ in the order topology such that $p \in U \subset B_\delta(p)$. 

  First, we claim that there is some $t\in I_{p,\infty} \cup \{\infty\}$ such that $I_{p,t}\subset B_\delta(p)$. Suppose that $I_{p,\infty}$ is nonempty; otherwise, we take $t=\infty$ and check that $I_{p,t} =\emptyset \subset B_\delta(p)$. By completeness, $I_{p,\infty}$ has a greatest lower bound $q$. If $q\ne p$, then $I_{p,q}=\emptyset$, so we suppose that $p=q$. By the above, $p$ is a limit point of $I_{p,\infty}$ in the subspace topology, i.e., for every $r>0$, the set $B_r(p)\cap I_{p,\infty}$ is nonempty.

  Let $C$ be as in Lemma~\ref{lem:compact-intersections} and let $t\in B_{C^{-1}\delta}(p)\cap I_{p,\infty}$. For all $v\in I_{p,t}$, we have $p\prec v\prec t$, so by Lemma~\ref{lem:compact-intersections}, $I_{p,t} \subset B_\delta(p)$, as desired.

  Likewise, there is some $s\in \{-\infty\} \cup I_{-\infty,p}$ such that $I_{s,p}\subset B_\delta(p)$. Then $s\prec p\prec t$, and
  $$I_{s,t} = I_{s,p}\cup \{p\}\cup I_{p,t}\subset B_\delta(\0).$$
  Thus the order topology is finer than the subspace topology, so the two topologies coincide.  
\end{proof}

We will use these lemmas to prove Proposition~\ref{prop:intervals}. 

\begin{proof}[Proof of Proposition~\ref{prop:intervals}]
  Let $E\subset \HH$ be a connected $\lambda$--vertical subset. We claim that $E$ is homeomorphic to an interval.
If we prove the assertion when $E$ is closed, then we are done. Indeed, if $E$ is not closed, then $E\subset\overline{E}$ is homeomorphic to a connected subset of an interval, and thus $E$ itself is homeomorphic to an interval, as desired.

Hence, we assume $E$ is closed. Since $E$ is connected, Lemma~\ref{lem:Properties}.(\ref{Item3}) implies that $E$ is totally ordered. Thus, by Lemma~\ref{lem:topologies}, the order on $E$ is complete and the order topology on $E$ agrees with the subset topology.

We claim that the order on $E$ is separable and dense. Since $E$ is separable in the subset topology, it is separable in the order topology. The density of the order follows from the connectedness of $E$. That is, if $p,q\in E$ and $p\prec q$, then $I_{p,q}$ is nonempty; otherwise $I_{-\infty,q}$ and $I_{p,\infty}$ would be disjoint open subsets of $E$ that cover $E$. 

The conclusion then follows from a classical result due to Cantor on complete, separable, dense orders. Let $E'$ be $E$ without its endpoints, if any, and let $\mathcal{D}$ be a countable dense subset of $E'$. By a theorem of Cantor, any such subset is order-isomorphic to $\mathbb Q\cap (0,1)$. The completeness of $E$ lets us extend this isomorphism to an isomorphism from $E'$ to the interval $(0,1)$. This isomorphism is a homeomorphism from $E'$ to an open interval, and we conclude the proof by adding the endpoints, if any.
\end{proof}

It remains to show that $E$ can be parameterized by a locally bi-H\"older map. We need the following lemma.

\begin{lemma}\label{lem:Biholder}
  There is a $0<\theta<1$ depending on $\lambda$ such that the following holds. 
  Let $E$ be a connected $\lambda$--vertical curve.
  Let $v, w\in E$ such that $v\prec w$ and let $0 < \epsilon < 1$. There is an $N\in \mathbb{N}$ depending on $\lambda$ and $\epsilon$ such that there is a sequence $v=q_0\prec q_1\prec \ldots\prec q_n=w$ with $n\leq N$ that satisfies
  \begin{equation}\label{eqn:biHolder-upper}
    d(q_i,q_{i+1})\leq \epsilon d(v,w)\qquad \text{for $i=0,\dots,n-1$}
  \end{equation}
  and
  \begin{equation}\label{eqn:biHolder-lower}
    d(q_i,q_j)\geq \theta \epsilon d(v,w)\qquad \text{for all $i,j$.}
  \end{equation}
\end{lemma}
\begin{proof}    
  Let us fix $v,w\in E$ with $v\prec w$. We first construct a sequence $v=q_0\prec q_1\prec \ldots \prec q_n=w$ such that
  \begin{equation}\label{eqn:biHolder-intermediate}
    \frac{\epsilon}{2} d(v,w) \le d(q_i,q_{i+1})\leq \epsilon d(v,w)\qquad \text{for $i=0,\dots,n-1$};
  \end{equation}
  then we will prove \eqref{eqn:biHolder-lower} and bound $n$.

  We proceed inductively. Let $q_0 = p$. Suppose that $i\ge 0$ and we have already constructed $q_i$ such that $v\preceq q_i\prec w$. If $d(q_i,w) \le \epsilon d(v,w)$, we let $q_{i+1}=w$ and terminate the process. Otherwise, $w$ is outside the ball $B(q_i,\frac{\epsilon}{2} d(v,w))$. The interval $\{x\in E: q_i\preceq x\preceq w\}$ is connected, so there is a point $q_{i+1}\in \partial B(q_i,\frac{\epsilon}{2} d(v,w))$ such that $q_i\prec q_{i+1}\preceq w$. This results in a sequence of $q_i$ that satisfies \eqref{eqn:biHolder-intermediate} and either terminates ($q_n=w$ for some $n$) or is infinite.

  We claim that $q_n=w$ for some $n$ and that we can bound $n$ by a function of $\lambda$.
  By Lemma \ref{lem:Properties}.(\ref{Item3}), $E$ is totally ordered, so $q_i\prec q_j$ for all $i<j$. Let $i<j$. Then $q_i\prec q_{i+1}\preceq q_j$ and $d(q_i,q_{i+1})\ge \frac{\epsilon}{2}$. By Lemma~\ref{lem:Properties}.(\ref{Item4}), there is a $c_\lambda>0$ such that 
  \begin{equation}\label{eqn:biHolder-constlower}
    d(q_i,q_j)\ge c_\lambda d(q_i,q_{i+1}) \ge c_\lambda\frac{\epsilon}{2} d(v,w),
  \end{equation}
  i.e., \eqref{eqn:biHolder-lower} holds.

  For every $q_i$, we have $v\prec q_i\prec w$, i.e., $q_i\in v\VCone^+_\lambda\cap w \VCone^-_\lambda$.
  By Lemma~\ref{lem:compact-intersections}, this intersection is compact. By \eqref{eqn:biHolder-constlower}, the balls of radius $c_\lambda\frac{\epsilon}{4}d(v,w)$ around the $q_i$'s are disjoint. If $C$ is as in Lemma~\ref{lem:compact-intersections}, all of these balls are contained in a ball of radius $(C+c_\lambda)d(v,w)$.
  Let
  $$
  N = (C+c_\lambda)^4\left(c_\lambda\frac{\epsilon}{4}\right)^{-4}.
  $$
  The ball of radius $(C+c_\lambda)d(v,w)$ contains at most $N$ disjoint balls of radius $c_\lambda\frac{\epsilon}{4}d(v,w)$, so $q_n=w$ for some $n\le N$. This proves the lemma.
\end{proof}

\begin{lemma}\label{lem:biHolderonCompact}
  Let $E$ be a connected $\lambda$--vertical curve. There are exponents $0<\alpha<\beta<1$ such that any compact $K\subset E$ is parametrized by a map $\gamma\from [0,1]\to K$ which is $(\alpha,\beta)$--bi-H\"older on compact sets, i.e., there is a $C>0$ such that for all $s,t\in [0,1]$,
  $$C^{-1}|s-t|^\beta < d(\gamma(s),\gamma(t)) < C |s-t|^\alpha.$$
\end{lemma}
\begin{proof}
    Let $\epsilon=\frac{1}{2}$, and let $\theta,N$ be as in Lemma~\ref{lem:Biholder}. Let $D$ be the set of fractions $kN^{-\ell}$ with $\ell,k$ integers such that $\ell\geq 0$, and $0\leq k\leq N^{\ell}$. 
    
    Let $v,w\in E$ and let $K=\{s\in E:v\preceq s\preceq w\}$. By iteratively applying Lemma \ref{lem:Biholder} we can construct an order preserving map $\gamma\from D\to K$ such that $\gamma(0)=v$, $\gamma(1)=w$, and for every $\ell\geq 0$ we have 
     \begin{equation}\label{eqn:Controlabove}
\max_{k=0,\ldots,N^\ell-1}d\left(\gamma\left(\frac{k}{N^\ell}\right),\gamma\left(\frac{k+1}{N^\ell}\right)\right)\leq 2^{-\ell} d(v,w), 
     \end{equation}
     \[
\min_{k=0,\ldots,N^\ell-1}d\left(\gamma\left(\frac{k}{N^\ell}\right),\gamma\left(\frac{k+1}{N^\ell}\right)\right)\geq \left(\frac{\theta}{2}\right)^\ell d(v,w).
\]
  The upper bound $d(\gamma(s),\gamma(t)) \lesssim |s-t|^{\log_N(2)}d(v,w)$ for all $s,t\in D$ is straightforward; see for instance \cite[Lemma 2]{LyonsVictoir}. For the lower bound, as in \cite[Page 90]{Kozhevnikov}, if $s,t\in D$, $s<t$, and $\ell$ is such that 
  $N^{-\ell} < \frac{1}{3}|s-t| \le N^{-\ell + 1},$
  then there are $k, \ell$ such that $s<k N^{-\ell} < (k+1) N^{-\ell} <t$. By applying \eqref{eq:coarse-monotone} twice, we find 
  \[
  \begin{aligned}
    d(s,t) &\ge c_\lambda^2 d\left(\gamma(k N^{-\ell}), \gamma((k+1) N^{-\ell})\right) \\ &\gtrsim \left(\frac{\theta}{2}\right)^\ell d(v,w) \gtrsim |s-t|^{-\log_N \frac{\theta}{2}}d(v,w).
    \end{aligned}
  \]
  Thus, $\gamma$ is bi-H\"older on $D$, and we can extend $\gamma$ to a bi-H\"older map on $[0,1]$ by continuity.

  Finally, we show that the image of $\gamma$ is $K$. Let $s\in K$. Since $K$ is totally ordered, for every $\ell$ we can find $k_\ell$ such that $\gamma(k_\ell N^{-\ell})\preceq s\prec \gamma((k_\ell+1)N^{-\ell})$. By Lemma \ref{lem:compact-intersections} and \eqref{eqn:Controlabove} we have that $d(\gamma(k_\ell N^{-\ell}),s)\to 0$ as $\ell\to \infty$. After passing to a convergent subsequence, suppose that $k_\ell N^{-\ell} \to t$ as $\ell\to \infty$; then $t\in[0,1]$ and $\gamma(t)=s$, as desired.
\end{proof}

Proposition~\ref{prop:Properties} then follows from Lemma~\ref{lem:Properties} and Lemma~\ref{lem:biHolderonCompact}.

\section{Vertical curves and intrinsic Lipschitz graphs}
\label{sec:containment}
In this section we study the relationship between vertical curves and intrinsic Lipschitz graphs and prove Lemma~\ref{lem:containment} and Theorem \ref{thm:intersections}. We first prove Lemma~\ref{lem:containment}.

\begin{proof}[Proof of Lemma~\ref{lem:containment}]
  First, if $(1 + 16 \lambda^2)^{-\frac{1}{4}}\le L < 1$, then
  \begin{equation}\label{eq:containment-condition-1}
    \VCone_\lambda \subset \Cone_{W,L} \qquad \text{for every vertical plane }W.
  \end{equation}
  Indeed, if $|z(p)| \ge \lambda|\pi(p)|^2$ for some $p\in\HH$, then
  $$d(\0,p)^4 = |\pi(p)|^4 + 16|z(p)|^2 \ge |\pi(p)|^4(1 + 16 \lambda^2) \ge |\pi(p)|^4 L^{-4},$$
  so for every vertical plane $W$,
  \[
  d(p,W) \leq |\pi(p)|\leq Ld_{\Kor}(\0,p).
  \]
  Thus, for every $\lambda>0$, there is an $L\in (0,1)$, and for every $L\in (0,1)$, there is a $\lambda>0$ such that \eqref{eq:containment-condition-1} holds. 
  
  Suppose that \eqref{eq:containment-condition-1} holds and that $E$ is a $\lambda$--vertical curve. Then, for any vertical plane $W$ and any $p\in E$, we have $E\subset p\VCone_\lambda \subset p\Cone_{W,L}$, so $E$ is an $L$--iLip graph over $W$; this, joined with Theorem \ref{thm:extension}, proves (1).

  To prove (2), we notice that if $0< L \le (1 + 16 \lambda^2)^{-\frac{1}{4}}$, then
  \begin{equation}\label{eq:containment-condition-2}
    \bigcap_W \Cone_{W,L} \subset \VCone_\lambda,
  \end{equation}
  where the intersection is taken over all vertical planes $W$. Indeed, let $p\in \bigcap_W \Cone_{W,L}$. Call $W_p:=|\pi(p)|^\perp$. Then $p\in \Cone_{W_p,L}$ from which we deduce
  \[
 d_{\Kor}(p,W_p)=|\pi(p)|\leq Ld_{\Kor}(\0,p)=L\sqrt[4]{|\pi(p)|^4+16|z(p)|^2}.
  \]
  From the latter inequality, we deduce as above that $\lambda|\pi(p)|^2\leq z|(p)|$, as desired.
  
  If $E\subset \HH$ and $E\subset p \Cone_{W,L}$ for all $W$, and all $p\in E$, then for all $p\in E$,
  $$
  E\subset\bigcap_W p\Cone_{W,L} \subset p\VCone_\lambda,
  $$
  and then $E$ is a $\lambda$--vertical curve. This proves (2).
\end{proof}

By Lemma~\ref{lem:containment} 
for any vertical curve $E$ and any vertical plane $W$, we can construct an entire iLip graph $\Gamma$ containing $E$. 

Second, when $L\in (0,1)$ is small, the intersection of two $L$--iLip graphs over different planes is a vertical curve.
\begin{lemma}\label{lem-connectedness}
    Let $W_1,W_2\subset\HH$ be two-dimensional vertical subgroups of $\HH$ such that $W_1\neq W_2$. There exists $L\in (0,1)$ such that if $\Gamma_1$ and $\Gamma_2$ are entire $L$--iLip graphs over $W_1$ and $W_2$ respectively, then $E=\Gamma_1\cap\Gamma_2$ is a nonempty connected vertical curve with no endpoints.
\end{lemma}
\begin{proof}
  Let $\lambda=\frac{1}{4}$ and let $B=B_1(\0)$ be the unit ball. Then
  $$\partial B \cap \bigcap_{L>0}(\Cone_{W_1,L} \cap \Cone_{W_2,L}) = \partial B \cap \langle Z\rangle$$
  lies in the interior of $\VCone_\lambda$. If $L\in (0,1)$ is small enough, then 
  $$\partial B \cap \Cone_{W_1,L} \cap \Cone_{W_2,L} \subset \VCone_\lambda,$$
  and since $\Cone_{W,L}$ and $\VCone_\lambda$ are scale-invariant, 
  \begin{equation}\label{eq:vertical-intersection}
  \Cone_{W_1,L} \cap \Cone_{W_2,L} \subset \VCone_\lambda.
  \end{equation}
  We choose $L$ small enough that \eqref{eq:vertical-intersection} holds and $\sin^{-1} L < \frac{1}{4}\angle(W_1,W_2)$.

  Let $\Gamma_1$ and $\Gamma_2$ be entire $L$--iLip graphs over $W_1$ and $W_2$ respectively and let $E=\Gamma_1\cap \Gamma_2$. If $p,q\in E$, then 
  $$
  q\in p\Cone_{W_1,L} \cap p \Cone_{W_2,L}\subset p\VCone_\lambda,
  $$
  so $E$ is $\lambda$--vertical.  
  We claim that $E$ is nonempty, connected, and has no endpoints.

  We first show that $E$ is nonempty. We will need some additional notation. The complement $\HH\setminus \Cone_{W_i,L}$ is the disjoint union of two open connected cones which we denote
  \[
     \HH\setminus \Cone_{W_i,L}:=\Cone_{V_i,L}^+\sqcup \Cone_{V_i,L}^-.
  \]
  Let $S^1$ be the horizontal unit circle $S^1=\{(x,y,0):x^2+y^2=0\}$. Let $V_i:=W_i^\perp$ and let $U_i\in V_i$ be a horizontal unit normal to $W_i$, chosen to point toward $\Cone_{V_i,L}^+$.

    \begin{figure}
  \begin{center}
  \begin{tikzpicture}
  \foreach \angle in {0, 60, 180, 240} {
    \draw[draw=none, fill=black!30] (0,0) -- (\angle-10:2.25) arc(\angle-10:\angle+10:2.25) -- cycle;
  }

  \draw[thick] (0, 0) circle (2cm);
  \foreach \angle/\label in {30/\phi_1, 120/\phi_2, 210/\phi_3, 300/\phi_4} {
    \node at (\angle:2) [circle, fill,inner sep=1.3pt,label={[label distance=-.1cm]\angle:$\label$}] {};
  }

  \draw[thick] (-2.25,0) -- (2.25,0); 
  \draw[thick] ({-2.25*cos(60)}, {-2.25 * sin(60)}) -- ({2.25*cos(60)}, {2.25 * sin(60)});
  \node at (2.25,0) [label={[label distance=-.2cm ]0:$W_1$}] {}; 
  \node at ({2.25*cos(60)}, {2.25 * sin(60)}) [label={[label distance=-.2cm]60:$W_2$}] {}; 

  \foreach \angle/\label in {75/I_1, 150/I_2, 255/I_3, 330/I_4} {
    \coordinate (p) at ({2*cos(\angle)},{2*sin(\angle)});
    \node at (p) [label={[label distance=-.2cm]{\angle+180}:$\label$}] {};
  }  
\end{tikzpicture}
  \end{center}
  \caption{\label{fig:cones}The intersection of $W_1$ and $W_2$ with the horizontal circle $S^1$. The shaded neighborhoods of $W_1$ and $W_2$ represent $\Cone_{W_1,L}$ and $\Cone_{W_2,L}$.  
  }
  \end{figure}
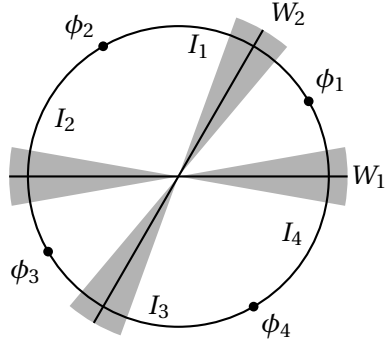
  
  For $\theta\in S^1$ and $t>0$, we have $d(\0,\theta^t) = t$ and $d(\theta^t, W_i) = t \sin \angle(\theta,W_i) = t |\cos\angle(\theta,U_i)|$, so 
  \begin{equation}\label{eq:plane-dist}
    d(\theta^t, W_i) - L d(\0,\theta^t) = t (\sin(\angle(\theta,W_i)) - L).
  \end{equation}
  For $p\in \HH$, let $D(p) = d(p, W_i) - L d(\0,p)$. Then $D$ is a $(1+L)$--Lipschitz function and $\Cone_{W_i,L} = D^{-1}((-\infty,0]))$; it follows that $\theta^t \in \Cone_{W_i,L}$ if and only if $\angle(\theta,W_i) \le \sin^{-1}(L)$ and that if $\sin(\angle(\theta,W_i))> L$, then 
  \begin{equation}\label{eq:cone-dist}
    d(\theta^t, \Cone_{W_i,L}) \ge \frac{t}{1+L}(\sin(\angle(\theta,W_i)) - L).
  \end{equation}
  
  By our choice of $L$, if $\angle(\theta,W_i) \ge \frac{1}{4}\angle(W_1,W_2)$, then $\theta^t \in \HH\setminus \Cone_{W_i,L}$, and whether $\theta^t\in \Cone_{V_i,L}^+$ or $\Cone_{V_i,L}^-$ depends on the sign of $\cos \angle(\theta, U_i)$.  
  
  The planes $W_1$ and $W_2$ divide $S_1$ into four arcs. Let $\phi_1, \phi_2, \phi_3, \phi_4 \in S^1$ be the midpoints of the four arcs in counterclockwise order, and let $I_j$ be the closed arc from $\phi_j$ to $\phi_{j+1}$ (taking $\phi_5=\phi_1$). Then $\angle(\phi_j,W_i) \ge \frac{1}{2}\angle(W_1,W_2)$ for all $i$ and $j$, so $\phi_j \in \Cone_{V_1,L}^\pm \cap \Cone_{V_2,L}^\pm$ for some choice of signs. After possibly swapping $\Cone_{V_i,L}^+$ and $\Cone_{V_i,L}^-$ and choosing the order of the $\phi_j$'s, we may suppose that $I_1\subset \Cone_{V_1,L}^+$, $I_2\subset \Cone_{V_2,L}^+$, $I_3\subset \Cone_{V_1,L}^-$, and $I_4\subset \Cone_{V_2,L}^-$, as in Figure~\ref{fig:cones}.

  Let $\Gamma_i^+$ and $\Gamma_i^-$ be the two components of $\HH\setminus \Gamma_i$, labeled so that $U_i$ points into $\Gamma_i^+$.
  Let $F=(f_1,f_2)\from \HH\to \R^2$, where the $f_i$'s are the signed distance functions
  $$f_i(v) = \begin{cases}
    d(v,\Gamma_i) & v \in \Gamma_i^+ \\
    0 & v \in \Gamma_i \\
    -d(v,\Gamma_i) & v \in \Gamma_i^-.
  \end{cases}$$
  We claim that if $t>0$ is sufficiently large, and $\gamma_t\from S^1\to \HH$, $\gamma_t(\theta)=\theta^t$, then the composition $F\circ \gamma_t$ has winding number $1$ around the origin.

  By \eqref{eq:cone-dist}, if $|\cos \angle(\theta, U_i)| > L$, then 
  $$d(\theta^t, \Cone_{W_i,L}) \ge \frac{t}{1+L} (|\cos \angle(\theta,U_i)| - L)  
  \ge \frac{t}{2} (|\cos \angle(\theta,U_i)| - L).$$
  Thus, if $g\in \Gamma_i$, $\theta\in S^1$, and $t>0$, then 
  \begin{equation}\label{eq:signed-on-graph}
    |f_i(g\theta^t)| \ge d(g\theta^t, g\Cone_{W_i,L}) \ge \frac{t}{2} (|\cos \angle(\theta,U_i)| - L),
  \end{equation}
  and the sign of $f_i(g\theta^t)$ agrees with the sign of $\cos \angle(\theta,U_i)$. 
  
  For any $h\in \HH$, we have
  $$g h g^{-1} h^{-1} = [g,h] = Z^{\omega(\pi(g),\pi(h))}$$
  and $|\omega(\pi(g),\pi(h))|\le \|g\|_K \|h\|_K$. In particular,
  $$d(h, g h g^{-1}) \le 2 \sqrt{\|g\|_K \|h\|_K}.$$
  Therefore,
  $$d(g \theta^{t}, \theta^t) = \|\theta^{-t} g^{-1} \theta^{t}\|_K  = \|[\theta^{-t},g^{-1}] \cdot g^{-1}\|_K \le 2\sqrt{|t|\|g\|_K} + d(\0, g) = O(\sqrt{t})$$
  as $t\to \infty$ with constant depending on $\|g\|_K$. If $|\cos \angle(\theta,U_i)|>L$, then the fact that $f_i$ is $1$--Lipschitz implies
  \begin{equation}\label{eq:signed-off-graph}
    |f_i(\theta^t)| \ge |f_i(g\theta^t)| - d(g\theta^t, \theta^t) \ge \frac{t}{2} (|\cos \angle(\theta,U_i)| - L) - O(\sqrt{t}),
  \end{equation}
  which is positive when $t$ is sufficiently large.

  A compactness argument combined with \eqref{eq:signed-off-graph} implies that if $t$ is sufficiently large, then
  $f_1(\gamma_t(I_1))>0$, $f_2(\gamma_t(I_2))>0$, $f_1(\gamma_t(I_3))<0$, and $f_2(\gamma_t(I_4))<0$. Therefore, $F(\gamma_t(S_1))\subset \R^2\setminus \0$ and $F\circ \gamma_t$ has winding number $1$ around $\0$. It follows that if $\Sigma$ is a surface in $\HH$ with boundary $\gamma_t$, then there is some $p\in \Sigma$ such that $F(p)=\0$ and thus $p\in \Gamma_1\cap \Gamma_2$. Thus $E=\Gamma_1\cap \Gamma_2$ is nonempty.

  Now, we claim that $E$ is connected and has no endpoints. Suppose that $\0\in E$. Then for any $\delta>0$, \eqref{eq:signed-on-graph} implies that $f_1(\gamma_\delta(I_1))>0$, $f_2(\gamma_\delta(I_2))>0$, $f_1(\gamma_\delta(I_3))<0$, and $f_2(\gamma_\delta(I_4))<0$, so any surface with boundary $\gamma_\delta$ contains an element of $E$. In particular, $\gamma_\delta$ divides the boundary $\partial B(\0,\delta)$ into two halves  $\Sigma^+=\{v\in \partial B(\0,\delta) : v\succ \0\}$ and $\Sigma^-=\{v\in \partial B(\0,\delta) : v\prec \0\}$, and each half contains an element of $E$. That is, there are points $q^\pm \in E$ such that $d(\0,q^\pm)=\delta$ and $q^-\prec \0 \prec q^+$. (Since $\lambda \ge \frac{1}{4}$, by Lemma \ref{lem:Properties}.(\ref{Item3}), $E$ is totally ordered.) By left-invariance, for any $p\in E$, there are $q^\pm \in E$ such that $d(p,q^\pm)=\delta$ and $q^-\prec p \prec q^+$.

  This property lets us describe the structure of $E$. By Lemma~\ref{lem:topologies}, the order on $E$ is complete and separable. We claim that it is dense. Suppose that $v, w\in E$ and $v\prec w$. Then $z(v^{-1}w)>0$. By Lemma~\ref{lem:Properties}.(\ref{Item4}), there is a $c_\lambda>0$ such that if $v\prec w \prec q$, then $z(v^{-1}q)\ge c_\lambda z(v^{-1}w)$ and thus
  $$d(v,q) \ge 2\sqrt{c_\lambda z(v^{-1}w)}.$$
  Let $\delta = \sqrt{c_\lambda z(v^{-1}w)}$. By the above, there is a $q\in E$ such that $v\prec q$ and $d(v,q)=\delta$. By our choice of $\delta$, we have $v\prec q \prec w$. Thus the order on $E$ is complete, separable, and dense. As in the proof of Proposition~\ref{prop:intervals}, $E$ is order-isomorphic to $(0,1)$ and thus it is homeomorphic to $(0,1)$, i.e., it is connected and does not contain any endpoints.
\end{proof}

\begin{proof}[Proof of Theorem \ref{thm:intersections}]
    Let $L'\in (0,1)$ be small enough such that Lemma~\ref{lem-connectedness} holds, i.e., the intersection of any $L'$--iLip graph over $W_1$ and any $L'$--iLip graph over $W_2$ is a nonempty connected vertical curve without endpoints. We may suppose that $L'$ is small enough that $\Cone_{W_1,L'}\cap\Cone_{W_2,L'}$ contains no horizontal vectors.

    By Lemma \ref{lem:containment}, there is a $\lambda_0>0$ such that any $\lambda_0$--vertical curve is contained in an entire $L'$--iLip graph over any two-dimensional vertical subgroup. We claim that the theorem holds for any $\lambda > \lambda_0$.

    Let $\lambda > \lambda_0$. 
    Let $S$ be the unit sphere in $\HH$. For any $L$, $\Cone_{W_1,L}\cap \Cone_{W_2,L}$ is a closed set, and as $L\to 0$, 
    $$S \cap \left(\Cone_{W_1,L} \cap \Cone_{W_2,L}\right)\to S\cap \langle Z\rangle = \{Z^{\pm\frac{1}{4}}\}.$$
    Since $S$ is compact and $\VCone_\lambda$ contains a neighborhood of $S\cap \langle Z\rangle$, there is some $L = L(\lambda)\in (0,L')$ such that
    $$S \cap \left(\Cone_{W_1,L} \cap \Cone_{W_2,L}\right)\subset \VCone_\lambda.$$
    Thus, by scaling invariance
\begin{equation}\label{ContainmentCones2}
\Cone_{W_1,L}\cap\Cone_{W_2,L}\subset  \VCone_\lambda.
    \end{equation}
    We claim that $L$ satisfies part (1) of the theorem. Let $\Gamma_1$ and $\Gamma_2$ be entire $L(\lambda)$-iLip graphs over $W_1$ and $W_2$ respectively. Since $L<L'$, by Lemma~\ref{lem-connectedness}, $E=\Gamma_1\cap \Gamma_2$ is a nonempty connected vertical curve with no endpoints. 
    By \eqref{ContainmentCones2}, $E$ is $\lambda$--vertical.

    Now, we claim that $L'$ satisfies part (2) of the theorem. Let $E$ be a $\lambda$--vertical curve. By Lemma \ref{lem:containment} we can find two entire $L'$--iLip graphs $\Gamma_1$ and $\Gamma_2$ over $W_1$ and $W_2$ such that $E\subset \Gamma_1\cap\Gamma_2$. Moreover, by Lemma~\ref{lem-connectedness}, $\Gamma_1\cap\Gamma_2$ is a connected vertical curve without endpoints. 
    
    If furthermore $E$ is closed, connected, nonempty, and has no endpoints, then $E$ is both closed and open in $\Gamma_1\cap \Gamma_2$. Since $\Gamma_1\cap \Gamma_2$ is connected, this implies $E=\Gamma_1\cap \Gamma_2$ as desired.
\end{proof}

\section{Vertical curves with fractional Hausdorff dimension}\label{Proof2.5}
In this section, we construct vertical curves with fractional Hausdorff dimension and prove Theorem~\ref{thm:HausDim2}.

We start by defining some basic concepts for smooth vertical curves, including the speed and the $\mathfrak{h}$--slope.
For a smooth curve $\alpha\from I\to \HH$, we let $\alpha'\from I \to \mathfrak{h}$ denote the velocity vector of $\alpha$, left-translated to the origin, i.e.,
\begin{equation}\label{eq:left-deriv}
  \alpha'(\tau) = \frac{\mathrm{d}}{\mathrm{d}t}[\alpha(\tau)^{-1}\alpha(t)]_{t=\tau}.
\end{equation}
\begin{lemma}\label{lem:measure-calculation}
  There exists a universal constant $c>0$ such that the following holds. Let $\alpha:[a,b]\to \HH$ be a smooth curve. Then $\cH^2(\alpha) = c\int_a^b |z(\alpha'(t))|\ud t$.
\end{lemma}
\begin{proof}
    The following proof is inspired by the computations in \cite[page 83]{Kozhevnikov}.
    For every $\tau\in (a,b)$ we have
\[
\begin{aligned}
\lim_{t\to 0}&\frac{\mathrm{d}_{\Kor}(\alpha(\tau+t),\alpha(\tau))^2}{t} \\ &=\lim_{t\to 0} \frac{\sqrt{\big|\pi(\alpha(\tau+t))-\pi(\alpha(\tau))\big|^4+16z(\alpha(\tau)^{-1}\alpha(\tau+t))^2}}{t} \\ 
&= 4|z(\alpha'(\tau))|.
\end{aligned}
\]
As a result, we have that 
\begin{equation}\label{eqn:Expansion}
\mathrm{d}_{\Kor}(\alpha(w),\alpha(v))^2 = 4|z(\alpha'(u))||v-w|+o(|v-u|+|w-u|), \quad \forall u,v,w\in (a,b),
\end{equation}
compare with \cite[Theorem 2]{Kirchheim}. Arguing verbatim as in \cite[Theorem 7, Lemma 4]{Kirchheim}, and using \eqref{eqn:Expansion}, we have that there exists a universal constant $c>0$ such that 
\[
\mathcal{H}^2(\alpha)=c\int_a^b |z(\alpha'(t))|\mathrm{d}t,
\]
as desired. 
\end{proof}
\begin{defn}\label{eqn:RescaledHausdorff}
    Let $c>0$ be the constant for which Lemma \ref{lem:measure-calculation} holds. We define a rescaled $2$-Hausdorff measure as follows:
    \[
    \hat{\cH^2}:=c^{-1}\cH^2.
    \]
\end{defn}
Consequently, if $\alpha$ is a smooth curve with $|z(\alpha'(t))|=1$ for all $t$, we say that it has \emph{unit vertical speed}. A $\lambda$--vertical curve may have $|z(\alpha'(t))|=0$ at some points, so not every $\lambda$--vertical curve has a unit-vertical-speed parameterization.

If $\gamma$ is a smooth curve and $z(\gamma'(t)) > 0$ for all $t$, we say that $\gamma$ is \emph{positively oriented}. Any such curve can be reparameterized with unit vertical speed. If $\gamma$ is positively oriented, we define its \emph{$\mathfrak{h}$--slope} at $t$ as $\frac{z(\gamma'(t))}{|\pi(\gamma'(t))|}$ when $\pi(\gamma'(t))\ne 0$ and as $+\infty$ if $\pi(\gamma'(t)) = 0$. For $m\ge 0$, we say that $\gamma$ has $\mathfrak{h}$--slope at least $m$ if $\gamma$ is positively oriented and 
$z(\gamma'(t)) \ge m |\pi(\gamma'(t))|$ for all $t$. The $\mathfrak{h}$--slope is independent of the parameterization of $\gamma$, but rescaling a curve multiplies its $\mathfrak{h}$--slope by a constant factor.
\begin{lemma}\label{lem:slope-scale}
  Let $m, r > 0$ and let $\gamma$ be a curve with $\mathfrak{h}$--slope at least $m$. Then $s_r\circ \gamma$ has $\mathfrak{h}$--slope at least $rm$.
\end{lemma}
\begin{proof}
  For all $t$,
  $$z((s_r\circ \gamma)'(t)) = z(s_r(\gamma'(t))) = r^2 z(\gamma'(t))$$
  and
  $$\pi((s_r\circ \gamma)'(t)) = \pi(s_r(\gamma'(t))) = r \pi(\gamma'(t)),$$
  so 
  $$z((s_r\circ \gamma)'(t)) = r^2 z(\gamma'(t)) \ge m r^2 |\pi(\gamma'(t))| = m r |\pi((s_r\circ \gamma)'(t))|,$$
  as desired.
\end{proof}

A curve with positive slope is locally vertical, but not necessarily globally vertical.
\begin{lemma}\label{lem:slope-vertical}
  Let $m>0$ and let $\alpha\from I\to \HH$ be a curve with $\mathfrak{h}$--slope at least $m$. Let $s,t\in \R$ and let
  $$\sigma = \left|\int_{s}^{t} z(\alpha'(\tau))\ud\tau\right|.$$
  Then $|\pi(\alpha^{-1}(s)\alpha(t))| \le \frac{\sigma}{m}$ and 
  \begin{equation}\label{eq:slope-vertical-diff}
    \left|z(\alpha^{-1}(s)\alpha(t)) -\sigma\right| \le \frac{\sigma^2}{m^2}.
  \end{equation}
  Consequently, if $\sigma < m^2$, then
  \begin{enumerate}
  \item $\alpha^{-1}(s)\alpha(t) \in \VCone_\lambda$, where 
  $$\lambda = \frac{m^2}{\sigma} - 1,$$
  \item $d(\alpha(s)^{-1}\alpha(t), Z^{\sigma}) \le \frac{5\sigma}{m}$, so 
  $$d(\alpha(s),\alpha(t)) \le 4\sqrt{\sigma} + \frac{5\sigma}{m}.$$
  \end{enumerate}
\end{lemma}
Some condition on $\sigma$ is necessary, since one can construct closed curves with positive $\mathfrak{h}$--slope. For example, $\gamma(t) = (\cos t, - \sin t, -t)$
is a horizontal curve, so $\lambda(t) = \gamma(t)Z^{t} = (\cos t, -\sin t, 0)$ is a closed curve with $\mathfrak{h}$--slope $1$.

\begin{proof}[Proof of Lemma~\ref{lem:slope-vertical}]
  After reparameterizing and translating, and possibly swapping $s$ and $t$, we may suppose that $\alpha$ is a unit-vertical-speed curve, $s=0\le t$, and $\alpha(s)=\0$, so that $\sigma=t-s=t$. A direct computation shows that 
  \[
\alpha'(t)=\left(\alpha_x'(t),\alpha_y'(t),\alpha_z'(t)-\frac{1}{2}\alpha_x(t)\alpha_y'(t)+\frac{1}{2}\alpha_x'(t)\alpha_y(t)\right).
  \]
  Since $\alpha$ is parametrized with unit vertical speed and has $\mathfrak{h}$--slope at least $m$ we have $z(\alpha'(t)) = 1$ and $|\pi(\alpha'(t))|\leq m^{-1}$. Integrating this inequality, we get the desired inequality $|\pi(\alpha(t))|\leq tm^{-1}$.
  
  By Cauchy--Schwarz,   
  \begin{equation}\label{eq:az-diff}
  \big|\alpha_z'(t)-1\big|=\frac{1}{2} \big|\alpha_x(t)\alpha_y'(t) -\alpha_x'(t)\alpha_y(t)\big| \le \frac{1}{2} |\pi(\alpha'(t))|\cdot |\pi(\alpha(t))| \le \frac{t}{2m^2}.
  \end{equation}
  Therefore, 
  \begin{equation}
    |\alpha_z(t)-t|=\left|\int_0^t \left(\alpha_z'(s)-1\right)\mathrm{d}s\right|\leq\int_0^t \frac{s}{2m^2}\ud s= \frac{t^2}{4m^2}\leq \frac{t^2}{m^2},
  \end{equation}
  which implies \eqref{eq:slope-vertical-diff}. The last two assertions of the lemma follow from direct calculation using \eqref{eq:def-koranyi} and \eqref{eqn:Vcone}.
\end{proof}

The goal of this section is to prove the following proposition, which immediately implies Theorem~\ref{thm:HausDim2}. 
\begin{prop}\label{prop:curves}
  Let $\hat{\cH^2}$ be the rescaled Hausdorff measure defined as in Definition \ref{eqn:RescaledHausdorff}. Let $\lambda,\lambda',\epsilon > 0$, with $\lambda'<\lambda$, and let $\alpha\from [0,\ell]\to \HH$ be a unit-vertical-speed $\lambda$--vertical curve. Suppose that $\alpha$ has $\mathfrak{h}$--slope at least $m$ for some $m>0$. There are smooth $\lambda'$--vertical curves $\gamma_1,\gamma_2 \from [0,\ell] \to \HH$ connecting $\alpha(0)$ to $\alpha(\ell)$ such that $\hat{\cH^2}(\gamma_1) < \epsilon$, $\hat{\cH^2}(\gamma_2) >\epsilon^{-1}$, and $d_{\mathrm{Haus}}(\alpha,\gamma_i) < \epsilon$ for each $i$, where $d_{\mathrm{Haus}}$ denotes the Hausdorff distance. 
  
  Furthermore, there are $\lambda'$--vertical curves $\tau_1, \tau_2 \from [0,1] \to \HH$ connecting $\alpha(0)$ to $\alpha(\ell)$ such that $d_{\mathrm{Haus}}(\alpha,\tau_i) < \epsilon$, and $\dim_H(\tau_1)<2<\dim_H(\tau_2)$.
\end{prop}

The curves we will construct  look like helixes at many different scales, as in Figure~\ref{fig:Hausdim2}. The idea of this construction is that if $\beta$ is a curve with $\mathfrak{h}$--slope at least $m$ and if $\rho>0$ is large, then, by Lemma~\ref{lem:slope-scale}, the rescaling $s_{m^{-1} \rho} \circ \beta$ has $\mathfrak{h}$--slope at least $\rho$. By Lemma~\ref{lem:slope-vertical}, the intersection of $s_{m^{-1} \rho} \circ \beta$ with any unit ball is close to a vertical segment. Equivalently, the intersection of $\beta$ with any ball of radius $m\rho^{-1}$ is close to a vertical segment. 

Thus, given a starting curve $\gamma$ with $\mathfrak{h}$--slope bounded away from zero, there is some $r_1$ so that the intersection of $\gamma$ with any ball of radius $r_1$ is close to a vertical segment. We can perturb $\gamma$ by replacing these vertical segments with helixes, and this replacement can increase or decrease the measure of $\gamma$. Furthermore, if the perturbations are small, the perturbed curve still has $\mathfrak{h}$--slope bounded away from zero, so we can repeat the process at a smaller scale $r_2$ and so on to construct curves with arbitrarily large or small Hausdorff measure; these curves and their limits will satisfy Proposition~\ref{prop:curves}.

We will need the following lemmas.
\begin{lemma}\label{lem:inductive-curves}
  Let $\hat{\cH^2}$ be the rescaled Hausdorff measure defined as in Definition \ref{eqn:RescaledHausdorff}. There are $\beta, r>0$ ($\beta=\frac{1}{400}$ and $r>1000$ suffice) with the following property.
  Let $0<\kappa<1$. For any $\phi=\pm 1$, any compact interval $I$, and any smooth curve $\gamma\from I \to \HH$ with $\mathfrak{h}$--slope at least $\kappa^{-1}$ and $\hat{\cH^2}(\gamma)\geq 1$, there is a smooth curve $\tilde{\gamma}\from I\to \HH$ with the same endpoints as $\gamma$ which satisfies the following properties:
  \begin{enumerate}
  \item $\tilde{\gamma}$ has $\mathfrak{h}$--slope at least $\kappa^{-1} r^{-1}$.
  \item For all $t$, $d(\gamma(t),\tilde{\gamma}(t))\le \kappa r^{-1}$.
  \item If $\phi = 1$, then for all $t$,
    $$(1+2\beta \kappa^2) z(\gamma'(t)) \le z(\tilde{\gamma}'(t)) \le (1+6\beta \kappa^2) z(\gamma'(t)).$$
    If $\phi = -1$, then for all $t$,
    $$(1-6\beta \kappa^2) z(\gamma'(t)) \le z(\tilde{\gamma}'(t)) \le (1-2\beta \kappa^2) z(\gamma'(t)).$$
  \end{enumerate}
\end{lemma}

We defer the proof of Lemma~\ref{lem:inductive-curves} to the end of the section.

\begin{lemma}\label{lem:leeway}
  Let $0 < \lambda' < \lambda$. There is a $\delta > 0$ such that if $v\in \VCone_\lambda$ and $g,g'\in B(\0,\delta \sqrt{|z(v)|})$, then
  $g v g' \in \VCone_{\lambda'}$.
\end{lemma}
\begin{proof}
  Let $K=\{w\in \VCone_\lambda : |z(w)|=1\}$. This set is a compact subset of the interior of $\VCone_{\lambda'}$, so there is a $\delta>0$ such that 
  $$B(\0,\delta) \cdot K \cdot B(\0,\delta) \subset \VCone_{\lambda'}.$$
  We claim that this $\delta$ satisfies the lemma. Let $v\in \VCone_\lambda$. The lemma is trivial when $v=\0$, so we may suppose that $|z(v)|\ne 0$. Let $s=s_{\sqrt{|z(v)|}}$ so that $|z(s(v))|=1$. By the scale-invariance of $\VCone_\lambda$,
  $$s^{-1}(g v g')= s^{-1}(g)s^{-1}(v)s^{-1}(g') \in B(\0,\delta) \cdot K\cdot B(\0,\delta),$$
  so $s^{-1}(g v g') \in \VCone_{\lambda'}$ and thus $g v g'\in \VCone_{\lambda'}$, as desired.
\end{proof}

Given these lemmas, we can prove Proposition~\ref{prop:curves}.
\begin{proof}[Proof of Proposition~\ref{prop:curves}]
  Let $c=c_\lambda>0$ be as in Lemma~\ref{lem:Properties}.(\ref{Item4}); we can take $c<\frac{1}{4}$. 
By Lemma~\ref{lem:leeway}, there is a $\delta>0$ such that if $v\in \VCone_\lambda$ and $g,g'\in B(\0,\delta\sqrt{|z(v)|})$, then $g v g' \in \VCone_{\lambda'}$.   Let 
$$
\kappa = \min \left\{(\lambda+1)^{-2},\frac{1}{10}, \frac{\delta\sqrt{c}}{4},\frac{\epsilon}{2}\right\}
$$

  After a rescaling and reparametrization we may suppose $\alpha$ is a positively oriented unit-vertical-speed curve, that $\ell = \hat{\cH^2}(\alpha) \ge 1$, and that $\alpha$ has $\mathfrak{h}$--slope at least $\kappa^{-1}$. We use Lemma~\ref{lem:inductive-curves} to construct a sequence of curves $\alpha_0,\alpha_1,\dots\from [0,\ell] \to \HH$ as follows. Let $\beta=\frac{1}{400}$ and $\rho=1000$ so that Lemma~\ref{lem:inductive-curves}. holds with $r=\rho$. Take $\phi=\pm1$.   
  Let $\alpha_0=\alpha$. 

  Suppose that we have constructed $\alpha_i$ such that $\alpha_i\from [0,\ell]\to \HH$ is a smooth curve with $\mathfrak{h}$--slope at least $m_i := \kappa^{-1} \rho^{-i}$ and $\ell_i:=\hat{\cH^2}(\alpha_i)\ge 2^{-i}$. This is already satisfied for $i=0$. Let $\gamma_i(t) = s_{\rho^{i}} \circ \alpha_i$; then by Lemma~\ref{lem:slope-scale}, $\gamma_i$ has $\mathfrak{h}$--slope at least $\kappa^{-1}$ and $\hat{\cH^2}(\gamma_i) \ge 1$, so we can apply Lemma~\ref{lem:inductive-curves} with $r=\rho$ to $\gamma_i$ to obtain a curve $\tilde{\gamma}_i$. Let
  $\alpha_{i+1}(t) = s_{\rho^{-i}}\circ \tilde{\gamma}_i.$

  Let $\ell_{i+1}:=\hat{\cH^2}(\alpha_{i+1})$. Then $\ell_{i+1}\ge \frac{1}{2} \ell_i\ge 2^{-i-1}$ and, by Lemma~\ref{lem:slope-scale}, $\alpha_{i+1}$ has $\mathfrak{h}$--slope at least $\kappa^{-1}\rho^{-i-1}$, so we can repeat this process inductively to construct $\alpha_i$ for all $i$. Furthermore, using Lemma~\ref{lem:inductive-curves} for all $i$:
  \begin{enumerate}
  \item $\alpha_i$ has $\mathfrak{h}$--slope at least $\kappa^{-1}\rho^{-i}$.
  \item For all $t$, $d(\alpha_i(t),\alpha_{i+1}(t)) \le \kappa \rho^{-i-1}$. Therefore, for all $j<i$,
    \begin{equation}\label{eq:alpha-diff}
      d(\alpha_j(t),\alpha_{i}(t)) \le \sum_{n=j}^{i-1} d(\alpha_n(t),\alpha_{n+1}(t)) \le 2\kappa \rho^{-j-1}.
    \end{equation}
  \item If $\phi = 1$, then for all $t$,
    \begin{equation}\label{eq:zalpha-plus}
      (1+2\beta \kappa^2) z(\alpha_i'(t)) \le z(\alpha_{i+1}'(t)) \le (1+6\beta \kappa^2) z(\alpha_i'(t)).
    \end{equation}
    If $\phi = - 1$, then for all $t$,
    \begin{equation}\label{eq:zalpha-minus}
      (1-6\beta \kappa^2) z(\alpha_i'(t)) \le z(\alpha_{i+1}'(t)) \le (1-2\beta \kappa^2) z(\alpha_i'(t)).
    \end{equation}
  \end{enumerate}

  Next, we introduce some notation and prove some metric bounds on the $\alpha_i$ in preparation to show that each curve $\alpha_i$ is $\lambda'$--vertical. For each $i$, let 
  \begin{equation}\label{eq:sigma-j}
    \Sigma_i(t) = \int_{0}^{t} z(\alpha_i'(\tau))\ud\tau,
  \end{equation}
  so that $\Sigma_i(\ell)=\ell_i$.
  Since $\alpha_i$ has positive slope, $\Sigma_i$ is monotone increasing. 
  
  Let $0\le s \le t \le \ell$ and for each $i$, let
  $$\sigma_i=\Sigma_i(t)-\Sigma_i(s)= \int_{s}^{t} z(\alpha_i'(\tau))\ud\tau$$
  and let $v_i=\alpha_i(s)^{-1}\alpha_i(t)$.
  
  We want to bound $v_i$; we start with the case that $s$ and $t$ are close together. Suppose that $i\ge 0$ and $\sigma_i\le \rho^{-2i}$.
  Recall that $\alpha_i$ has $\mathfrak{h}$--slope at least $m_i=\kappa^{-1}\rho^{-i}$. By Lemma~\ref{lem:slope-vertical}, $v_i\in \VCone_L$, where
  $$
  L = \frac{m_i^2}{\rho^{-2i}} - 1 = \kappa^{-2} - 1 >\lambda$$
  and
  $$
  |z(v_i) - \sigma_i| \le \frac{1}{m_i^2}\sigma_i^2 \le \frac{\rho^{-2i}}{\kappa^{-2}\rho^{-2i}} \sigma_i\le \frac{1}{2} \sigma_i.
  $$
  Therefore, 
  \begin{equation}\label{eq:slope-bound}
    \frac{\sigma_i}{2}\le |z(v_i)| \le \frac{3\sigma_i}{2} \text{ and } v_i \in \VCone_\lambda \qquad \text{if $\sigma_i\le \rho^{-2i}$.}
  \end{equation}

  By choosing $j$ carefully, we can use \eqref{eq:slope-bound} to get upper and lower bounds on $v_j$. Suppose that $|s-t| \le 1$ and let $j$ be the largest integer such that $\sigma_j \le \rho^{-2j}$. By the maximality of $j$, $\sigma_{j+1} \ge \rho^{-2j-2}$.
  Regardless of $\phi$, \eqref{eq:zalpha-plus} and \eqref{eq:zalpha-minus} imply $z(\alpha_j'(\tau))\ge \frac{1}{2} z(\alpha_{j+1}'(\tau)),$
  so $\sigma_{j}\ge \frac{1}{2} \rho^{-2j-2}$. Combined with \eqref{eq:slope-bound}, this implies
  \begin{equation}\label{eq:z-upperlower}
  v_j\in \VCone_\lambda,\qquad  \frac{1}{2} \rho^{-2j-2}\le \sigma_j \le \rho^{-2j},\qquad
    \frac{\sigma_j}{2}\le z(v_j) \le\frac{3\sigma_j}{2}.
  \end{equation}

  We use \eqref{eq:z-upperlower} and Lemma~\ref{lem:leeway} to prove that the $\alpha_i$ are $\lambda'$--vertical. 

  Let $i\ge 0$, let $0\le s\le t\le \ell$, and let $v_i$ and $\sigma_i$ be as above. We claim that $v_i\in \VCone_{\lambda'}$; we consider three cases. First, if $\sigma_i \le \rho^{-2i}$, then $v_i\in \VCone_{\lambda}$ by \eqref{eq:slope-bound}. 
  
  Second, suppose that $\sigma_i > \rho^{-2i}$ but $|s-t|<1$. Let $j$ be the largest integer such that $\sigma_j \le \rho^{-2j}$; note that $j<i$. Then $v_j\in \VCone_\lambda$ by \eqref{eq:slope-bound}; we will use Lemma~\ref{lem:leeway} to show that $v_i\in \VCone_\lambda'$.
  
  By \eqref{eq:alpha-diff}, there are $g,g'\in B(\0,2\kappa\rho^{-j-1})$ such that $\alpha_i(s) = \alpha_j(s) g$ and $\alpha_i(t) = \alpha_j(t) g'$, so
  $$v_i = \alpha_i(s)^{-1} \alpha_i(t) = g^{-1} \alpha_j(s)^{-1}\alpha_j(t)g' = g^{-1}v_j g'.$$
  By \eqref{eq:z-upperlower}, $z(v_j)\ge \frac{1}{4} \rho^{-2j-2} \ge \frac{c}{2}\rho^{-2j-2}$, and by our choice of $\kappa$, 
  $$2\kappa \rho^{-j-1} \le \delta \sqrt{\frac{c}{2}}\rho^{-j-1} \le \delta \sqrt{|z(v_j)|}.$$
  Lemma~\ref{lem:leeway} then implies that $v_i\in \VCone_{\lambda'}$.
  
  Finally, suppose that $|s-t|\ge 1$. Then $v_0\in \VCone_\lambda$ by the verticality of $\alpha_0$. Since $\alpha_0$ is a unit-vertical-speed curve, we have $\sigma_0 = t - s \ge 1$, i.e., $s+1\le t$. Since $\alpha_0$ has $\mathfrak{h}$--slope at least $\kappa^{-1}$ and $\kappa^{-1}\ge 10$, Lemma~\ref{lem:slope-vertical} implies that
  $$z(\alpha_0(s)^{-1}\alpha_0(s+1)) \ge 1 - \frac{1}{100} \ge \frac{1}{2}.$$
  Since $\alpha_0$ is positively oriented, $\alpha_0(s)\preceq \alpha_0(s+1)\preceq \alpha_0(t)$, so by \eqref{eq:coarse-monotone}, $z(v_0) \ge c z(\alpha_0(s)^{-1}\alpha_0(s+1)) \geq \frac{c}{2}$. We conclude as above; by \eqref{eq:alpha-diff}, there are $g, g'\in B(\0,2\kappa\rho^{-1})$ such that $v_i = g^{-1}v_0g'$, so by Lemma~\ref{lem:leeway}, $v_i\in \VCone_{\lambda'}$. These three cases are the only possibilities, so $\alpha_i$ is $\lambda'$--vertical for all $i$.

  When $\phi=1$, we have $\ell_{i+1}\ge (1+2\beta\kappa^2) \ell_i$, so $\ell_i\to \infty$ as $i\to \infty$; when $\phi = -1$, $\ell_{i+1}\le (1-2\beta\kappa^2) \ell_i$, so $\ell_i\to 0$ as $i\to \infty$. Thus, when $\phi=1$, this construction produces a family of $\lambda'$--vertical curves with arbitrarily large $\hat{\cH^2}$--measure, and when $\phi=-1$, it produces curves with arbitrarily small $\hat{\cH^2}$--measure. Taking also into account that $\kappa\leq \epsilon/2$, and \eqref{eq:alpha-diff}, this proves the first part of the proposition.

  It remains to construct vertical curves with Hausdorff dimension greater or less than $2$. For all $t$, let $\zeta(t) = \lim_{i\to\infty} \alpha_i(t)$. This limit exists thanks to \eqref{eq:alpha-diff}. Since $\VCone_{\lambda'}$ is closed, $\zeta$ is $\lambda'$--vertical.

  Let
  $$m_\phi =
  \begin{cases}
    1-6\beta \kappa^2 & \phi = -1 \\
    1+2\beta \kappa^2 & \phi = 1,
  \end{cases}$$
  $$M_\phi =
  \begin{cases}
    1-2\beta \kappa^2 & \phi = -1 \\
    1+6\beta \kappa^2 & \phi = 1,
  \end{cases}$$
  so that we can combine \eqref{eq:zalpha-plus} and \eqref{eq:zalpha-minus} into
  \begin{equation}\label{eq:mphi}
    m_\phi z(\alpha_i'(t)) \le z(\alpha_{i+1}'(t)) \le M_\phi z(\alpha_i'(t)).
  \end{equation}

  Let
  $$c_\phi = \frac{\log \rho}{\log (\rho^2M_\phi)} = \frac{1}{2+ \log_\rho M_\phi}$$
  $$C_\phi = \frac{\log \rho}{\log (\rho^2m_\phi)} = \frac{1}{2+ \log_\rho m_\phi}.$$
  We claim that for all $s, t\in [0,\ell]$ such that $s<t$ and $|s-t|\le 1$, we have
  \begin{equation}\label{eq:zeta-holder}
    |s-t|^{C_\phi} \lesssim_\lambda d(\zeta(s),\zeta(t)) \lesssim_\lambda |s-t|^{c_\phi}.
  \end{equation}

  As above, for all $i$, let $\sigma_i =\Sigma_i(t)-\Sigma_i(s)$ and $v_i=\alpha_i(s)^{-1}\alpha_i(t)$, and let $j$ be the largest integer such that $\sigma_j \le \rho^{-2j}$. By \eqref{eq:z-upperlower},  $\frac{1}{4} \rho^{2j-2}\le z(v_j)\le \frac{3}{2}\rho^{-2j}$ and $v_j\in \VCone_\lambda$, so
  $$\frac{1}{2}\rho^{-j-1} \le \sqrt{|z(v_j)|} \le d(\alpha_j(s),\alpha_j(t)) \lesssim_\lambda \rho^{-j}.$$
  By \eqref{eq:alpha-diff}, 
$d(\zeta(\tau), \alpha_j(\tau)) \le 2\kappa \rho^{-j-1}$ for all $\tau$, so $\frac{1}{4} \rho^{-j-1}\le d(\zeta(s),\zeta(t)) \lesssim_\lambda \rho^{-j}$, i.e., $d(\zeta(s),\zeta(t)) \approx_\lambda \rho^{-j}$.

  By \eqref{eq:mphi},
  $$m_\phi^{j}|s-t|\le \sigma_j \le M_\phi^{j}|s-t|.$$
  Since $\sigma_j\approx \rho^{-2j}$,
  $$(\rho^{2}M_\phi)^{-j} \lesssim |s-t| \lesssim (\rho^{2}m_\phi)^{-j}.$$
  Therefore,
  $$|s-t|^{c_\phi} \gtrsim (\rho^{2}M_\phi)^{-jc_\phi} = \rho^{-j} \approx_\lambda d(\zeta(s),\zeta(t))$$
  and
  $$|s-t|^{C_\phi} \lesssim (\rho^{2}m_\phi)^{-jC_\phi} = \rho^{-j} \approx_\lambda d(\zeta(s),\zeta(t)),$$
  which proves \eqref{eq:zeta-holder}. Thus, $\zeta$ is a locally bi-Hölder map, and
  $$\dim_H(\zeta(I)) \in [C_\phi^{-1}, c_\phi^{-1}] = [2 + \log_\rho m_\phi, 2+ \log_\rho M_\phi].$$
  If $\phi=1$, this implies $\dim_H(\zeta(I))>2$; if $\phi=-1$, this implies $\dim_H(\zeta(I))<2$, as desired. 
\end{proof}

\begin{proof}[Proof of Lemma \ref{lem:inductive-curves}]   
    Let us fix $\beta=\frac{1}{400}$, and take $r>1000$.
    
    We first reparametrize $\gamma\from I\to\HH$ with unit vertical speed. Let $\ell=\hat{\cH^2}(\gamma)\geq 1$. Suppose that $I=[a,b]$ and let $\sigma\from I \to \R$,
    $$\sigma(t) = \int_a^t z(\gamma'(\tau))\ud \tau.$$
    Then $\sigma'(t)>0$ for all $t$ and $\sigma(b)=\ell$, so we define $\lambda\from [0,\ell]\to \HH$, $\lambda(t) = \gamma(\sigma^{-1}(t))$. This is a smooth, unit vertical speed curve.
    
    We apply a perturbation with wavelength approximately $2\pi r^{-2}$ to $\lambda$.
    Let $\alpha = \big\lceil\frac{r^2 \ell}{2\pi}\big\rceil$ and $\xi=\frac{2\pi \alpha}{\ell}$. Since $r^2\ell > 1000^2$, we have $\frac{r^2 \ell}{2\pi} \le \alpha\le\frac{r^2 \ell}{\pi}$ and $r^2 \le \xi \le 2r^2$.
    
    Let $v\from [0,\ell] \to \HH$,
    $$v(t) = \frac{\kappa}{8r} (\cos(\xi t)X -\phi\sin(\xi t)Y).$$
    Let $w(s) = v(0)^{-1} v(s)$ and let
    \begin{equation}\label{eqn:tildelambda}
    \tilde \lambda(s):=\lambda(s)\cdot w(s).
    \end{equation}
    Then $w(0)=w(\ell)=\0$, so the endpoints of $\tilde\lambda$ and $\lambda$ coincide. 
    
    We claim that $\tilde{\lambda}$ satisfies versions of properties (1)--(3), i.e., $\tilde{\lambda}$ has $\mathfrak{h}$--slope at least $\kappa^{-1}r^{-1}$, $d(\lambda(t),\tilde{\lambda}(t))\le \kappa r^{-1}$, and $z(\tilde{\lambda}'(t))$ is between $1 + 2\phi \beta \kappa^2$ and $1 + 6 \phi\beta \kappa^2$ for all $t$.
    
    First, note that $d(\0,v(t)) = \frac{\kappa}{8}r^{-1}$ for all $t$, so $d(\lambda(t),\tilde{\lambda}(t))=d(\0,w(t))\leq \kappa r^{-1}$. This proves the second property.

    To prove the first and third properties, we bound $|\pi(\tilde{\lambda}'(t))|$ and $z(\tilde{\lambda}'(t))$.
    By \eqref{eq:left-deriv} and \eqref{eqn:tildelambda}
  \begin{align*}
    \tilde{\lambda}'(\tau) 
    & =
    \frac{\mathrm{d}}{\mathrm{d}t} [w(\tau)^{-1}\lambda(\tau)^{-1} \lambda(t) w(t) ]_{t=\tau} \\
    & =
    \frac{\mathrm{d}}{\mathrm{d}t} [w(\tau)^{-1}\lambda(\tau)^{-1}\lambda(t) w(\tau) \cdot w(\tau)^{-1} w(t) ]_{t=\tau} \\
    & = \Adj_{w(\tau)^{-1}}(\lambda'(\tau)) + w'(\tau),
  \end{align*}
  where, for $v\in \mathfrak h$ and $g\in\HH$, we denote $\Adj_{g}(v) := \frac{\mathrm{d}}{\mathrm{d}t}[g \exp(tv) g^{-1}]_{t=0}\in\mathfrak{h}$, or, since we use exponential coordinates to identify $\mathfrak{h}$ and $\HH$, $\Adj_{g}(v) := \frac{\mathrm{d}}{\mathrm{d}t}[g v^t g^{-1}]_{t=0}$.

  By \eqref{eq:group-law}, for $v\in \mathfrak h$ and $g\in\HH$,
  $$(g v^t) g^{-1} = \left(g + tv + \tfrac{1}{2}[g,tv]\right) - g + \tfrac{1}{2}\left[g + tv + \tfrac{1}{2}[g,tv],-g\right] = t v + [g,tv].$$
  Letting $\omega$ be the area form on $\R^2$ as in Section~\ref{sec:H1},
  $$\Adj_g(v) = v + [g,v] = v + \omega(\pi(g),\pi(v)) Z.$$
  In particular, $\Adj_{g}(v) - v\in \langle Z\rangle$ and
  $$
  |\Adj_{g}(v) - v| \leq |\pi(g)| |\pi(v)|.
  $$
  Then, 
  \begin{equation}\label{eq:pi-diff}
    \pi(\tilde{\lambda}'(t)) = \pi\left(\Adj_{w(t)^{-1}}(\lambda'(t)) + w'(t)\right) = \pi(\lambda'(t)) + \pi(w'(t)),
  \end{equation}
  and
  \begin{multline}\label{eq:z-diff}
    \left|z(\tilde{\lambda}'(t)) - \left(z(\lambda'(t)) + z(w'(t))\right)\right| = \left| z\left(\Adj_{w(t)^{-1}}\left(\lambda'(t)\right) + w'(t) - \lambda'(t) - w'(t)\right)\right| \\
    = \left|z\left(\Adj_{w(t)^{-1}}\left(\lambda'(t)\right) - \lambda'(t)\right)\right| \le \big|\pi(w(t))\big| \big|\pi(\lambda'(t))\big|.
  \end{multline}
  
  Furthermore,
  \begin{equation}
  w'(t) = v'(t) = \frac{\kappa\xi r^{-1}}{8} v'(\xi t)+\frac{\kappa^{2}r^{-2}\xi}{128}\phi Z,
  \end{equation}
  so, using that $\ell\ge \alpha\pi r^{-2}$,
  \begin{equation}\label{eq:pi-wt}
    |\pi(w'(t))| =\frac{\kappa\xi r^{-1}}{8}=\frac{2\kappa\alpha\pi r^{-1}}{8\ell} \leq \frac{2\kappa\alpha\pi r^{-1}}{8\alpha\pi r^{-2}}\leq \frac{1}{4}\kappa r.
  \end{equation}

  Since $r^2 \le \xi \le 2r^2$, we have 
  \begin{equation}\label{eq:z-wtphi}
  \frac{\kappa^2}{128}\leq |z(w'(t))| =\frac{\kappa^2 r^{-2}\xi}{128}  \leq \frac{\kappa^2}{64},
  \end{equation}
  and $z(w'(t))$ has the same sign as $\phi$. 

  This lets us bound $\tilde{\lambda}'(t)$. 
  Since $\lambda$ is unit-vertical-speed and has $\mathfrak{h}$--slope at least $\kappa^{-1}$, we have $z(\lambda'(t))=1$ and $|\pi(\lambda'(t))| \le \kappa$ for all $t$. By \eqref{eq:pi-diff} and \eqref{eq:pi-wt}, and since $r>4$,
  \begin{equation}\label{est-pi-gtild}
  |\pi(\tilde{\lambda}'(t))| \le |\pi(\lambda'(t))| + |\pi(w'(t))|\leq \kappa+\frac{1}{4}\kappa r \leq \frac{1}{2}\kappa r. 
  \end{equation}
  By \eqref{eq:z-diff} we have
  \begin{equation}\label{eqn:Est-gamma'}
  \begin{aligned}
  \left|z(\tilde{\lambda}'(t)) - 1 -z(w'(t))\right| &\le |\pi(w(t))|\cdot |\pi(\lambda'(t))| \\&\le \kappa r^{-1}\cdot \kappa= \kappa^2 r^{-1}.
  \end{aligned}
  \end{equation}
  
  Since $\beta=\frac{1}{400}$ and $r>1000$, if $\phi=-1$, then \eqref{eqn:Est-gamma'} and \eqref{eq:z-wtphi} imply
  \begin{equation*}
  1 - 6\beta \kappa^2 \le 1 - \frac{\kappa^2}{64}-\frac{\kappa^2}{r} \leq z(\tilde\lambda '(t)) \leq 1 - \frac{\kappa^2}{128} +\frac{\kappa^2}{r} \le 1 - 2\beta \kappa^2.
  \end{equation*}
  If $\phi = 1$, then 
  \begin{equation*}
  1 + 2\beta \kappa^2 \le 1 + \frac{\kappa^2}{128}-\frac{\kappa^2}{r} \leq z(\tilde\lambda '(t)) \leq 1 + \frac{\kappa^2}{64} +\frac{\kappa^2}{r} \le 1 + 6\beta \kappa^2.
  \end{equation*}
  In either case, $z(\tilde{\lambda}'(t))$ is between $1 + 2 \phi \beta \kappa^2$ and $1 + 6 \phi \beta \kappa^2$. In particular, $z(\tilde{\lambda}'(t)) \ge \frac{1}{2}$, so by \eqref{est-pi-gtild} the $\mathfrak{h}$-slope of $\tilde\lambda$ is at least $\frac{1}{2}\left(\frac{1}{2}\kappa r\right)^{-1}=\kappa^{-1}r^{1}$, as desired.

  Finally, for all $t\in I$, let $\tilde{\gamma}(t) = \tilde{\lambda}(\sigma(t))$. Since $\tilde{\gamma}$ is a reparametrization of $\tilde{\lambda}$, it has the same endpoints as $\gamma$ and satisfies properties (1)--(3).
\end{proof}

\subsection{Comparison with previous works}\label{sec:Previous} In \cite[Chapters 5 \& 6]{Kozhevnikov}, and \cite{KozhevnikovArticle}, A. Kozhevnikov studied properties of a class of vertical curves in the $(2n+1)$-dimensional Heisenberg group $\mathbb H^{2n+1}$. See also \cite{LeonardiMagnani, MagnaniStepanovTrevisan} for related results.

In \cite{Kozhevnikov}, Kozhevnikov considered curves of the following form. Let $F\from  \mathbb H^{2n+1}\to\mathbb R^{2n}$ of regularity $C^1_{\mathrm{H}}$ such that $F(\0)=\0$, $\Gamma=F^{-1}(\0)$, and $\mathrm{D}_{v}F$ is surjective for all $v\in \Gamma$. Here $\mathrm{D}_v F$ denotes the Pansu differential of $F$ at $v$. Then $F^{-1}(0)$ is a $1$--manifold:
\begin{thm}[Theorem 5.3.7 \cite{Kozhevnikov}]\label{corvertical}
    Let $F\in C^1_{\mathrm H}(\mathbb H^{2n+1},\mathbb R^{2n})$ be such that $F(\0) = \0$ and $\mathrm{D}_\0 F$ is surjective. Then there is a neighborhood $U$ of 0 such that $U\cap F^{-1}(\0)$ is a simple curve.
\end{thm}
Furthermore, Kozhevnikov shows that if $U$ is sufficiently small, then $\Gamma = U\cap F^{-1}(\0)$ satisfies some additional properties, including the fact that \eqref{eq:def-order} is a total ordering, that $\Gamma$ is vanishing Reifenberg flat, and that $\Gamma$ has a bi-Hölder parameterization; he called sets with these properties \emph{vertical curves}.

To avoid confusion with the $\lambda$--vertical curves in this paper, we call sets of the form $F^{-1}(\0)$ \emph{$C^1_{\mathrm{H}}$--fibers}. If a set of the form $U\cap F^{-1}(\0)$ satisfies Kozhevnikov's additional conditions, we call it a \emph{regular vertical curve}.

We stress that the notion of $\lambda$--vertical curve is \textbf{locally} weaker than that of a regular vertical curve. If $E\subset \HH$ is a regular vertical curve, then, by the fact that $E$ is vanishing Reifenberg flat (see Remark 5.3.2 \cite{Kozhevnikov}), for every $p\in E$ there are $\delta,\lambda>0$ such that $B_\delta(p)\cap E$ is a $\lambda$--vertical curve. Globally, however, $C^1_{\mathrm{H}}$--fibers need not be $\lambda$--vertical curves and vice versa.

We further point out that Theorem \ref{corvertical} is in analogy with the more general Proposition \ref{prop:intervals}. Moreover, in the proof of Theorem \ref{corvertical} the author proves that any regular vertical curve is a bi-H\"older curve, compare with \cite[Theorem 5.2.14 \& Theorem 5.3.5]{Kozhevnikov}, and with our Lemma \ref{lem:biHolderonCompact}. Further parametrization results related to Theorem \ref{corvertical}, and Proposition \ref{prop:intervals} are in \cite[Theorem 1.1]{LeonardiMagnani}, and \cite[Theorem 5.6]{MagnaniStepanovTrevisan}.
\smallskip

Kozhevnikov also studied properties related to the Hausdorff dimension of vertical curves. In particular he proved Theorem \ref{thm:KozHaus} below, and constructed pathological examples (Proposition \ref{PathologicalExamples} below) that inspired our Proposition \ref{prop:curves}. 

\begin{thm}[Proposition 5.3.10 \& Corollary 5.4.8 \& Corollary 5.4.16 in \cite{Kozhevnikov}]\label{thm:KozHaus}
    The Hausdorff dimension of a regular vertical curve is 2.
\end{thm}
\begin{prop}[Example 5.6.16, Example 5.6.17, Example 5.6.19 in \cite{Kozhevnikov}]\label{PathologicalExamples}
    There exist regular vertical curves $\Gamma_1,\Gamma_2$ such that $\mathcal{H}^2(\Gamma_1)=0$, and $\mathcal{H}^2(\Gamma_2)=\infty$. Moreover, there exists a $2$-Ahlfors regular vertical curve $\Gamma_3$ such that: for every $a\in \Gamma_3$
    \[
    0<\liminf_{r\to 0}\frac{\mathcal{H}^2(\Gamma_3\cap B_r(a))}{2r^2}<\limsup_{r\to 0}\frac{\mathcal{H}^2(\Gamma_3\cap B_r(a))}{2r^2}=1.
    \]
\end{prop}

Let us remark that as a result of the proofs of \cite[Lemma 5.4.3 \& Corollary 5.4.6 \& Corollary 5.4.8]{Kozhevnikov}, or refining our Lemma \ref{lem:biHolderonCompact}, it can be inferred that $|\dim_H(E)-2|\leq o_{\lambda\to\infty}(1)$ for every $\lambda$--vertical curve $E\subset \HH$.  But, in contrast with the result proved by Kozhevnikov in the Theorem \ref{thm:KozHaus}, an arbitrary $\lambda$--vertical curve might not have Hausdorff dimension 2, as we showed in Proposition \ref{prop:curves}.

%% file: contactdiffeo.tex
\section{Fibers of contact diffeomorphisms}\label{sec:FiberContact}

Let $B$ be the closed unit ball in $\HH$. In this section we will prove Theorem~\ref{thm:multiscale-maps} by constructing contact diffeomorphisms $\beta\from B\to B$ which are arbitrarily close to the identity map and have the property that the average of the $\mathcal{H}^2$-measure of the fibers of $\pi\circ\beta$ is arbitrarily small. Here $\pi\from \HH\to\mathbb R^2$ is the projection onto the horizontal plane.

We will need a few preliminaries in order to prove Theorem \ref{thm:multiscale-maps}. For a smooth map $f\from \HH \to \R^2$, the differential $Df_v$ sends $\mathfrak{h}$ to $\R^2$. We define the \emph{horizontal Jacobian} $J^H_f(v) := \det(Df_v|_{\langle X_v,Y_v\rangle})$ to be the determinant of $Df_v|_{\langle X_v, Y_v\rangle}$. The following co-area formula was proved for the first time in \cite{MagnaniCoarea}. For generalizations, see \cite{KarmanovaVodopyanov, Kozhevnikov, MagnaniStepanovTrevisan}.
\begin{thm}\label{thm:Coarea}
  There is a $c>0$ such that if $f\from \HH \to \R^2$ is a smooth map and $U\subset \HH$ is a measurable set, then
  \begin{equation}\label{eq:coarea}
    \int_{\R^2} \mathcal{S}^2(f^{-1}(v)\cap U) \ud v = c \int_U \left|J^H_f(g)\right| \ud g.
  \end{equation}
\end{thm}

Recall from Section~\ref{sec:Contact} that a contact diffeomorphism is a diffeomorphism $\beta \from \HH \to \HH$ that preserves the horizontal distribution, i.e., for any $v\in \HH$, $D_v\beta$ sends $\langle X_v, Y_v\rangle$ to $\langle X_{\beta(v)},Y_{\beta(v)}\rangle$, where $X$, $Y$, and $Z$ denote the standard left-invariant vector fields as in \eqref{eq:left-inv}. We define the \emph{horizontal Jacobian} of $\beta$ as $J^H_\beta(v) := \det(D_v\beta|_{\langle X,Y\rangle})$. By the coarea formula, to prove Theorem~\ref{thm:multiscale-maps}.(\ref{it:small-fibers}), it suffices to find $\beta$ such that
\begin{equation}\label{eq:jacobian-bound}
  \int_{B}\left| J^H_{\pi\circ \beta}(v)\right|\ud v = \int_{B}\left| J^H_\beta(v)\right|\ud v < \epsilon.
\end{equation}
We will show that one can reduce the left-hand side of \eqref{eq:jacobian-bound} by composing $\beta$ with a contact diffeomorphism supported on a small ball, then use the Vitali Covering Theorem to prove Theorem~\ref{thm:multiscale-maps}. We recall that, for a diffeomorphism $f\from A\to A$, $\supp f = \closure \{x\in A\mid f(x)\ne x\}$; we say that $f$ is supported on a ball $B$ if $\supp f \subset \inter(B)$. 

\begin{lemma}\label{lem:reduction}
  There is a $\kappa >0$ with the following property. Let $\beta\from \HH \to \HH$ be a contact map and let $p\in \HH$ be a point such that $J^H_\beta(p)\ne 0$. There is an $r_0 > 0$ such that for any $0<r<r_0$, there is a contact diffeomorphism $\alpha\from \HH \to \HH$ supported on $B(p,r)$ such that
  \begin{equation}\label{eq:reduction}
    \int_{B(p,r)}\left| J^H_{\beta\circ \alpha}(v)\right|\ud v \le (1 - \kappa) \int_{B(p,r)}\left| J^H_\beta(v)\right|\ud v.
  \end{equation}
\end{lemma}
\begin{proof}  
  We first construct a contact diffeomorphism $\alpha_0\from \HH \to \HH$ supported on $B:=B(\0,1)$ such that 
  $$\int_{B}\left| J^H_{\alpha_0}(v)\right|\ud v < \vol(B).$$
  In fact, it suffices to choose any $\alpha_0$ which is not volume-preserving.

  Recall from Section~\ref{sec:Contact} that any smooth function $\psi\in C^\infty(\HH)$ corresponds to a vector field 
  $$
  V_\psi := Y[\psi] X - X[\psi] Y + \psi Z,
  $$
  such that $V_\psi$ generates a contact flow, i.e., a flow $\Psi_\psi^t\from \HH \to \HH$ such that each map $\Psi_\psi^t$ is a contact diffeomorphism. Fix some nonzero $\psi\in C^\infty_c(B)$, so that for all $t$, $\Psi_\psi^t$ is a contact diffeomorphism supported on $B$.
  By \eqref{eqn:JacobianDet} we can choose $\psi\in C^\infty_c(B)$ and $t$ small enough so that $\Psi_\psi^t$ is not volume-preserving. Let $\alpha_0 := \Psi_\psi^t$.

  Since $\alpha_0$ is a contact diffeomorphism, we have $D_v\alpha_0(\langle X_v,Y_v\rangle) = \langle X_{\alpha_0(v)},Y_{\alpha_0(v)}\rangle$ for all $v$; let
  $$D_v^H\alpha_0 = \begin{pmatrix} a(v) & c(v) \\
    b(v) & d(v) 
  \end{pmatrix} \in M^{2\times 2}(\R).$$
  be the corresponding matrix. Then
  $$(\alpha_0)_*(Z) = [(\alpha_0)_*(X), (\alpha_0)_*(Y)] = [a X + b Y, c X + d Y].$$
  Since $[fV,gW] = fg[V,W] + fV[g] W - gW[f] V$, this has $Z$--component
  $$\mathrm{d}z\left((\alpha_0)_*(Z)\right) = ad - bc = J^H_{\alpha_0}(v),$$
  and
  \begin{equation}\label{eq:squares}
    J_{\alpha_0}(v) = J^H_{\alpha_0}(v)\cdot (\alpha_0)_*(Z)_Z = J^H_{\alpha_0}(v)^2.
  \end{equation}
  Furthermore, since $\alpha_0$ is a diffeomorphism, we have $J_{\alpha_0}(v)>0$ for all $v$ and thus $J^H_{\alpha_0}(v)>0$ for all $v$.

  Since $\alpha_0$ fixes $\partial B$ pointwise, the coarea formula for smooth maps implies that
  $$\int_B J_{\alpha_0}(v) \ud v = \int_B \ud v = \vol(B).$$
  Then, by Jensen's inequality,
  $$\frac{1}{\vol(B)} \int_B J^H_{\alpha_0}(v) \ud v \stackrel{\eqref{eq:squares}}{=} \frac{1}{\vol(B)} \int_B \sqrt{J_{\alpha_0}(v)} \ud v \le \sqrt{\frac{1}{\vol(B)} \int_B J_{\alpha_0}(v) \ud v} = 1.$$
  Since $\alpha_0$ is not volume-preserving, $J_{\alpha_0}$ is not constant on $B$, and this inequality is strict. That is,
  $$\int_B J^H_{\alpha_0}(v) \ud v < \vol(B).$$

  Let
  $$\kappa = \frac{1}{2}\left(1 - \frac{\int_B J^H_{\alpha_0}(v) \ud v}{\vol(B)}\right),$$
  so that 
  $$\int_{B}J^H_{\alpha_0}(v)\ud v = (1-2\kappa)\vol(B).$$
  
  Let $\beta\from \HH \to \HH$ be a contact map and let $p\in \HH$ be a point such that $J^H_\beta(p)\ne 0$. For $r>0$, let
  $$\alpha_{p,r}(v) = p\cdot (s_r\circ \alpha_0\circ s_r^{-1})(p^{-1} v),$$
  so that $\alpha_{p,r}$ is a contact diffeomorphism supported on $B(p,r)$. Then, making the substitution $w = s_r^{-1}(p^{-1} v)$, 
  \begin{multline*}
    \int_{B(p,r)} J^H_{\beta\circ \alpha_{p,r}}(v) \ud v = \int_{B(p,r)} J^H_{\beta}(\alpha_{p,r}(v)) J^H_{\alpha_{p,r}}(v)  \ud v \\
    = r^4 \int_{B} J^H_{\beta}(p\cdot s_r(\alpha_0(w))) J^H_{\alpha_{0}}(w) \ud w.
  \end{multline*}
  By the continuity of $J^H_{\beta}$, 
  $$\lim_{r\to 0} r^{-4} \int_{B(p,r)} J^H_{\beta\circ \alpha_{p,r}}(v) \ud v = \int_{B} J^H_{\beta}(p) J^H_{\alpha_{0}}(w) \ud w = (1-2\kappa) J^H_{\beta}(p)\vol(B)$$
  and
  $$\lim_{r\to 0} r^{-4} \int_{B(p,r)} J^H_{\beta}(v) \ud v = J^H_{\beta}(p)\vol(B).$$
  Therefore, when $r$ is sufficiently small, 
  $$\int_{B(p,r)} J^H_{\beta\circ \alpha_{p,r}}(v) \ud v \le (1 - \kappa) \int_{B(p,r)} J^H_\beta(v),$$
  as desired.
\end{proof}

\begin{proof}[Proof of Theorem~\ref{thm:multiscale-maps}]
  Since $\mathcal{H}^2\leq \mathcal{S}^2$, see \cite[Section 2.10.2]{FedererBook}, it is enough to prove the theorem with $\mathcal{S}^2$ instead of $\mathcal{H}^2$. We construct $\beta$ inductively. Let $\kappa$ be as in Lemma~\ref{lem:reduction} and let $B\subset \HH$ be the closed unit ball. Let $\beta_0 = \id_{\HH}$. We will define a sequence $A_0,A_1,A_2,\dots\from \HH \to\HH$ of contact diffeomorphisms such that if $\beta_{i+1} = \beta_i \circ A_i$ for all $i$ and 
  $$F_i := \int_{\R^2} \mathcal{S}^2((\pi\circ \beta_i)^{-1}(v) \cap B) \ud v = c \int_B \left|J^H_{\beta_i}(v)\right| \ud v,$$
  then
  \begin{enumerate}
  \item There is a finite union $U_i$ of disjoint balls, each of radius at most $\epsilon 2^{-i-1}$ such that $\supp A_i \subset \inter(U_i)$ and $U_i\subset \inter(B)$.
  \item $F_{i+1} \le (1 - \frac{\kappa}{2}) F_i$ for all $i$.
  \end{enumerate}

  For each $i$, we use the Vitali Covering Theorem to construct $U_i$. For each $p\in \HH$, let $r_0(p)$ be as in Lemma~\ref{lem:reduction} applied to $\beta_i$. For $p\in \inter(B)$, let
  $$0 < r(p) < \min\{d(p,\partial B), r_0(p), \epsilon 2^{-i-1}\}.$$
  The balls of the form $B(p,r(p))$ form a Vitali covering of $\inter(B)$, so by the Vitali Covering Theorem, there is a sequence of points $p_1,p_2,\dots\in \HH$ such that the balls $B(p_j,r(p_j))$ are all disjoint and $B\setminus \bigcup_j B(p_j,r(p_j))$ has measure zero. Let $U_{i,N} = \bigcup_{j=1}^N B(p_j,r(p_j))$. Then 
  $$\lim_{N\to \infty} \int_{U_{i,N}} \left|J^H_{\beta_i}(v)\right| \ud v = F_i;$$
  we choose $N$ large enough that $\int_{U_{i,N}} \left|J^H_{\beta_i}(v)\right| \ud v \ge \frac{F_i}{2}$ and let $U_i = U_{i,N}$. 

  For $j=1,\dots, N$, let $\alpha_j$ be as in Lemma~\ref{lem:reduction}, so that $\alpha_j$ is supported on $B(p_j,r(p_j))$ and
  \begin{equation}\label{eq:apply-reduce}
    \int_{B(p_j,r_j)}\left| J^H_{\beta_i\circ \alpha_j}(v)\right|\ud v \le (1-\kappa) \int_{B(p_j,r_j)}\left| J^H_{\beta_i}(v)\right|\ud v.
  \end{equation}
  Let $A_i\from \HH\to \HH$,
  $$A_i(v) =\begin{cases}
    \alpha_j(v) & v\in B(p_j,r(p_j)) \text{ for some } 1\le j \le n \\
    v & \text{otherwise}.
  \end{cases}$$
  Then $U_i$ and $A_i$ satisfy condition (1). Let $\beta_{i+1} = \beta_i \circ A_i$. By \eqref{eq:apply-reduce}, for $j=1,\dots, N$,
  $$\int_{B(p_j,r_j)}\left| J^H_{\beta_i}(v)\right| - \left| J^H_{\beta_{i+1}}(v)\right|\ud v \ge \kappa \int_{B(p_j,r_j)}\left| J^H_{\beta_i}(v)\right|\ud v.$$
  It follows that
  \begin{multline*}
    F_i - F_{i+1} = \int_B \left|J^H_{\beta_i}(v)\right| - \left|J^H_{\beta_{i+1}}(v)\right| \ud v \\ = \sum_{j=1}^N  \int_{B(p_j,r(p_j))} \left|J^H_{\beta_i}(v)\right| - \left|J^H_{\beta_{i+1}}(v)\right| \ud v \ge \kappa \int_{U_i} \left|J^H_{\beta_i}(v)\right|\ud v 
    \ge \frac{\kappa}{2} F_i.
  \end{multline*}
  and thus $F_{i+1} \le (1 - \frac{\kappa}{2}) F_i$ as desired.
\end{proof}

%% file: ref.bib
@article {KoranyiReimann,
    AUTHOR = {Kor\'anyi, A. and Reimann, H. M.},
     TITLE = {Foundations for the theory of quasiconformal mappings on the
              {H}eisenberg group},
   JOURNAL = {Adv. Math.},
  FJOURNAL = {Advances in Mathematics},
    VOLUME = {111},
      YEAR = {1995},
    NUMBER = {1},
     PAGES = {1--87},
      ISSN = {0001-8708,1090-2082},
   MRCLASS = {30C65 (31C15 43A80)},
  MRNUMBER = {1317384},
MRREVIEWER = {Gerald\ B.\ Folland},
       DOI = {10.1006/aima.1995.1017},
       URL = {https://doi.org/10.1006/aima.1995.1017},
}

@article {LyonsVictoir,
    AUTHOR = {Lyons, T. and Victoir, N.},
     TITLE = {An extension theorem to rough paths},
   JOURNAL = {Ann. Inst. H. Poincar\'e{} C Anal. Non Lin\'eaire},
  FJOURNAL = {Annales de l'Institut Henri Poincar\'e{} C. Analyse Non
              Lin\'eaire},
    VOLUME = {24},
      YEAR = {2007},
    NUMBER = {5},
     PAGES = {835--847},
      ISSN = {0294-1449,1873-1430},
   MRCLASS = {60H10 (60H07)},
  MRNUMBER = {2348055},
MRREVIEWER = {Antoine\ J.\ Lejay},
       DOI = {10.1016/j.anihpc.2006.07.004},
       URL = {https://doi.org/10.1016/j.anihpc.2006.07.004},
}

@article {AMiLip,
    AUTHOR = {Antonelli, G. and Merlo, A.},
     TITLE = {Intrinsically {L}ipschitz functions with normal target in
              {C}arnot groups},
   JOURNAL = {Ann. Fenn. Math.},
  FJOURNAL = {Annales Fennici Mathematici},
    VOLUME = {46},
      YEAR = {2021},
    NUMBER = {1},
     PAGES = {571--579},
      ISSN = {2737-0690,2737-114X},
   MRCLASS = {53C17 (22E25 26A16 28A75 49Q15)},
  MRNUMBER = {4277829},
MRREVIEWER = {Scott\ Robert\ Zimmerman},
       DOI = {10.5186/aasfm.2021.4638},
       URL = {https://doi.org/10.5186/aasfm.2021.4638},
}

@article {FSSCILip,
    AUTHOR = {Franchi, B. and Serapioni, R. and Serra Cassano,
              F.},
     TITLE = {Intrinsic {L}ipschitz graphs in {H}eisenberg groups},
   JOURNAL = {J. Nonlinear Convex Anal.},
  FJOURNAL = {Journal of Nonlinear and Convex Analysis. An International
              Journal},
    VOLUME = {7},
      YEAR = {2006},
    NUMBER = {3},
     PAGES = {423--441},
      ISSN = {1345-4773,1880-5221},
   MRCLASS = {58C20 (22E30)},
  MRNUMBER = {2287539},
MRREVIEWER = {Thierry\ Coulhon},
}

@article {FSIlip,
    AUTHOR = {Franchi, B. and Serapioni, R. P.},
     TITLE = {Intrinsic {L}ipschitz graphs within {C}arnot groups},
   JOURNAL = {J. Geom. Anal.},
  FJOURNAL = {Journal of Geometric Analysis},
    VOLUME = {26},
      YEAR = {2016},
    NUMBER = {3},
     PAGES = {1946--1994},
      ISSN = {1050-6926,1559-002X},
   MRCLASS = {49Q15 (22E25 53C17 58C20 58C35)},
  MRNUMBER = {3511465},
MRREVIEWER = {Jingzhi\ Tie},
       DOI = {10.1007/s12220-015-9615-5},
       URL = {https://doi.org/10.1007/s12220-015-9615-5},
}

@article {FSSCiLip2,
    AUTHOR = {Franchi, B. and Serapioni, R. and Serra Cassano,
              F.},
     TITLE = {Differentiability of intrinsic {L}ipschitz functions within
              {H}eisenberg groups},
   JOURNAL = {J. Geom. Anal.},
  FJOURNAL = {Journal of Geometric Analysis},
    VOLUME = {21},
      YEAR = {2011},
    NUMBER = {4},
     PAGES = {1044--1084},
      ISSN = {1050-6926,1559-002X},
   MRCLASS = {22E30 (58C20)},
  MRNUMBER = {2836591},
MRREVIEWER = {Davide\ Vittone},
       DOI = {10.1007/s12220-010-9178-4},
       URL = {https://doi.org/10.1007/s12220-010-9178-4},
}

@article {VittoneiLip,
    AUTHOR = {Vittone, D.},
     TITLE = {Lipschitz graphs and currents in {H}eisenberg groups},
   JOURNAL = {Forum Math. Sigma},
  FJOURNAL = {Forum of Mathematics. Sigma},
    VOLUME = {10},
      YEAR = {2022},
     PAGES = {Paper No. e6, 104},
      ISSN = {2050-5094},
   MRCLASS = {49Q15 (22E30 26A16 53C17)},
  MRNUMBER = {4377000},
MRREVIEWER = {Yongsheng\ Zhang},
       DOI = {10.1017/fms.2021.84},
       URL = {https://doi.org/10.1017/fms.2021.84},
}

@article {NY-VvsH,
    AUTHOR = {Naor, A. and Young, R.},
     TITLE = {Vertical perimeter versus horizontal perimeter},
   JOURNAL = {Ann. of Math. (2)},
  FJOURNAL = {Annals of Mathematics. Second Series},
    VOLUME = {188},
      YEAR = {2018},
    NUMBER = {1},
     PAGES = {171--279},
      ISSN = {0003-486X,1939-8980},
   MRCLASS = {46B85 (20F65)},
  MRNUMBER = {3815462},
MRREVIEWER = {Damian\ Sawicki},
       DOI = {10.4007/annals.2018.188.1.4},
       URL = {https://doi.org/10.4007/annals.2018.188.1.4},
}

@article{Kozhevnikov,
      title={{P}ropriétés métriques des ensembles de niveau des applications différentiables
		sur les groupes de {C}arnot}, 
      author={Kozhevnikov, A.},
      year={2015},
      journal = {Géométrie métrique [math.MG]. Université Paris Sud - Paris XI, (2015)}
}

@book {FedererBook,
    AUTHOR = {Federer, H.},
     TITLE = {Geometric measure theory},
    SERIES = {Die Grundlehren der mathematischen Wissenschaften},
    VOLUME = {Band 153},
 PUBLISHER = {Springer-Verlag New York, Inc., New York},
      YEAR = {1969},
     PAGES = {xiv+676},
   MRCLASS = {28.80 (26.00)},
  MRNUMBER = {257325},
MRREVIEWER = {J.\ E.\ Brothers},
}

@article {Kirchheim,
    AUTHOR = {Kirchheim, B.},
     TITLE = {Rectifiable metric spaces: local structure and regularity of
              the {H}ausdorff measure},
   JOURNAL = {Proc. Amer. Math. Soc.},
  FJOURNAL = {Proceedings of the American Mathematical Society},
    VOLUME = {121},
      YEAR = {1994},
    NUMBER = {1},
     PAGES = {113--123},
      ISSN = {0002-9939,1088-6826},
   MRCLASS = {28A78},
  MRNUMBER = {1189747},
MRREVIEWER = {G.\ Freilich},
       DOI = {10.2307/2160371},
       URL = {https://doi.org/10.2307/2160371},
}

@article{RigotiLip, 
    author={Rigot, S.},
    title={Quantitative notions of rectifiability in the Heisenberg groups},
    journal={arXiv:1904.06904.}
}

@article {KarmanovaVodopyanov,
    AUTHOR = {Karmanova, M. and Vodopyanov, S.},
     TITLE = {A coarea formula for smooth contact mappings of
              {C}arnot-{C}arath\'eodory spaces},
   JOURNAL = {Acta Appl. Math.},
  FJOURNAL = {Acta Applicandae Mathematicae},
    VOLUME = {128},
      YEAR = {2013},
     PAGES = {67--111},
      ISSN = {0167-8019,1572-9036},
   MRCLASS = {53C17 (51F99)},
  MRNUMBER = {3125636},
MRREVIEWER = {Davide\ Vittone},
       DOI = {10.1007/s10440-013-9822-7},
       URL = {https://doi.org/10.1007/s10440-013-9822-7},
}

@article{KozhevnikovArticle,
      title={Roughness of level sets of differentiable maps on {H}eisenberg group}, 
      author={Kozhevnikov, A.},
      year={2011},
      journal = {arXiv:1110.3634}
}

@article {LeonardiMagnani,
    AUTHOR = {Leonardi, G. P. and Magnani, V.},
     TITLE = {Intersections of intrinsic submanifolds in the {H}eisenberg
              group},
   JOURNAL = {J. Math. Anal. Appl.},
  FJOURNAL = {Journal of Mathematical Analysis and Applications},
    VOLUME = {378},
      YEAR = {2011},
    NUMBER = {1},
     PAGES = {98--108},
      ISSN = {0022-247X,1096-0813},
   MRCLASS = {53C17 (49Q20 53C40)},
  MRNUMBER = {2772447},
       DOI = {10.1016/j.jmaa.2010.12.052},
       URL = {https://doi.org/10.1016/j.jmaa.2010.12.052},
}

@article {MagnaniCoarea,
    AUTHOR = {Magnani, V.},
     TITLE = {Note on coarea formulae in the {H}eisenberg group},
   JOURNAL = {Publ. Mat.},
  FJOURNAL = {Publicacions Matem\`atiques},
    VOLUME = {48},
      YEAR = {2004},
    NUMBER = {2},
     PAGES = {409--422},
      ISSN = {0214-1493,2014-4350},
   MRCLASS = {28A75 (28A78 53C17)},
  MRNUMBER = {2091013},
MRREVIEWER = {Vasily\ A.\ Chernecky},
       DOI = {10.5565/PUBLMAT\_48204\_07},
       URL = {https://doi.org/10.5565/PUBLMAT_48204_07},
}

@article {MagnaniStepanovTrevisan,
    AUTHOR = {Magnani, V. and Stepanov, E. and Trevisan, D.},
     TITLE = {A rough calculus approach to level sets in the {H}eisenberg
              group},
   JOURNAL = {J. Lond. Math. Soc. (2)},
  FJOURNAL = {Journal of the London Mathematical Society. Second Series},
    VOLUME = {97},
      YEAR = {2018},
    NUMBER = {3},
     PAGES = {495--522},
      ISSN = {0024-6107,1469-7750},
   MRCLASS = {53C17 (28A75)},
  MRNUMBER = {3816397},
MRREVIEWER = {Davide\ Vittone},
       DOI = {10.1112/jlms.12115},
       URL = {https://doi.org/10.1112/jlms.12115},
}
